\documentclass[10pt]{article}

\usepackage{fullpage}
\usepackage{epsfig}
\usepackage{psfrag}
\usepackage{amssymb}
\usepackage{supertabular}
\newcommand{\old}[1]{{}}

\DeclareMathSymbol{\plusminus}{\mathbin}{symbols}{"06}
\renewcommand{\pm}{{{\plusminus}\mbox{\,}}}
\newcommand{\PM}{{\plusminus}}
\newcommand{\p}{{{+}\mbox{\,}}}
\newcommand{\m}{{{-}\mbox{\,}}}
\newcommand{\s}{{\ast}}
\renewcommand{\d}{{\cdot}}
\newcommand{\dd}{{:}}
\newcommand{\x}{{\times}}
\renewcommand{\b}{\bar}
\renewcommand{\o}{{\circ}}

\newcommand{\h}{\frac{1}{2}}
\newcommand{\q}{\frac{1}{4}}
\newcommand{\Q}{\frac{3}{4}}
\renewcommand{\t}{\frac{1}{3}}
\newcommand{\T}{\frac{2}{3}}
\newcommand{\e}{\frac{1}{8}}
\newcommand{\E}{\frac{3}{8}}

\newcommand{\up}[1]{\mbox{}^{\,#1}}

\newcommand{\<}{\langle}
\renewcommand{\>}{\rangle}

\newcommand{\wrapper}[1]{\multicolumn{2}{l}{\left.\begin{array}{ll}#1\end{array}\right .  }}

\newcommand{\QED}{\nobreak\hfill\hbox{$\square$}}

\newcommand{\subheader}[4]
{
	\hline\hline
	\multicolumn{4}{|l|}{Plane group: #1}\\
	\multicolumn{4}{|l|}{Relations #2 :}\\
	\multicolumn{4}{|p{6cm}|}{#3}\\
	\hline
	Fibrifold & Couplings for & Point & Intern.\\
	name & #4 & group & no.\\
	\hline
}
\newcommand{\n}{\,} 
\newcommand{\xo}{\hspace{1.3pt} 0}
\newcommand{\xh}{\frac{1}{2}}
\newcommand{\xt}{\frac{1}{3}}
\newcommand{\xT}{\frac{2}{3}}
\newcommand{\xq}{\frac{1}{4}}
\newcommand{\xQ}{\frac{3}{4}}
\newcommand{\xs}{\frac{1}{6}}

\newcommand{\xR}{\, 0\m}
\newcommand{\id}{($I$)}
\title{On Three-Dimensional Space Groups}

\author{
John H.~Conway\\
{\em Department of Mathematics}\\
{\em Princeton University, Princeton NJ 08544-1000 USA}\\
{\em e-mail: {\tt conway@math.princeton.edu}}\\
\and
Olaf Delgado Friedrichs\\
{\em Department of Mathematics}\\
{\em Bielefeld University, D-33501 Bielefeld}\\
{\em e-mail: {\tt delgado@mathematik.uni-bielefeld.de}}\\
\and
Daniel H.~Huson
\thanks{Current address: Celera Genomics, 45 West Gude Drive,
Rockville MD 20850 USA}
\\
{\em Applied and Computational Mathematics}\\
{\em Princeton University, Princeton NJ 08544-1000 USA}\\
{\em e-mail: {\tt huson@member.ams.org}}\\
\and
William P.~Thurston\\
{\em Department of Mathematics}\\
{\em University of California at Davis}\\
{\em e-mail: {\tt wpt@math.ucdavis.edu}}\\
}

\begin{document}

\maketitle

\begin{abstract}
A entirely new and independent enumeration of the crystallographic space
groups is given, based on obtaining the groups as fibrations over the plane
crystallographic groups, when this is possible.
For the 35 ``irreducible'' groups for which it is not, 
an independent method is used that has the advantage of
elucidating their subgroup relationships.
Each space group is given a short ``fibrifold name'' which,
much like the orbifold names for two-dimensional groups,
while being only specified up to isotopy, contains enough information
to allow the construction of the group from the name.
\end{abstract}

\section{Introduction}

There are 219 three-dimensional crystallographic space groups
(or 230 if we distinguish between mirror images).
They were independently enumerated in the 1890's by
W.~Barlow in England, E.S.~Federov in Russia and A.~Sch\"onfliess
in Germany. The groups are comprehensively described in
the International Tables for Crystallography \cite{Hahn83}.
For a brief definition, see Appendix~I.

Traditionally the enumeration depends on classifying lattices
into 14 Bravais types, distinguished by the symmetries that can be added,
and then adjoining such symmetries in all possible ways.
The details are complicated, because there are many cases to consider.

Here we present an entirely new and independent enumeration, based
on obtaining the groups as fibrations over the plane crystallographic
groups, when this is possible.
For the 35 ``irreducible'' groups for which it is not, 
we use an independent method that has the advantage of
elucidating their subgroup relationships.
We describe this first.

\section{The 35 Irreducible Groups}

A group is {\em reducible} or {\em irreducible} according as there is or is
not a direction that it preserves up to sign.

\subsection{Irreducible Groups}

We shall use the fact that any irreducible group $G$ has
elements of order 3, which generate what we call its {\em odd subgroup}
($3$ being the only possible odd order greater than $1$).
The odd subgroup $T$ of $G$ is obviously normal and so
$G$ lies between $T$ and its normalizer $N(T)$.

This is an extremely powerful remark, since it turns out that there are
only two possibilities $T_1$ and $T_2$ for the odd subgroup,
and $N(T_1)/T_1$ and $N(T_2)/T_2$ are finite groups of order
$16$ and $8$. This reduces the enumeration of irreducible space groups
to the trivial enumeration of subgroups of these two finite groups (up to
conjugacy).

The facts we have assumed will be proved in Appendix~II.

\subsection{The 27 ``Full Groups''}

\begin{figure}[c]
\hfil\epsfig{file=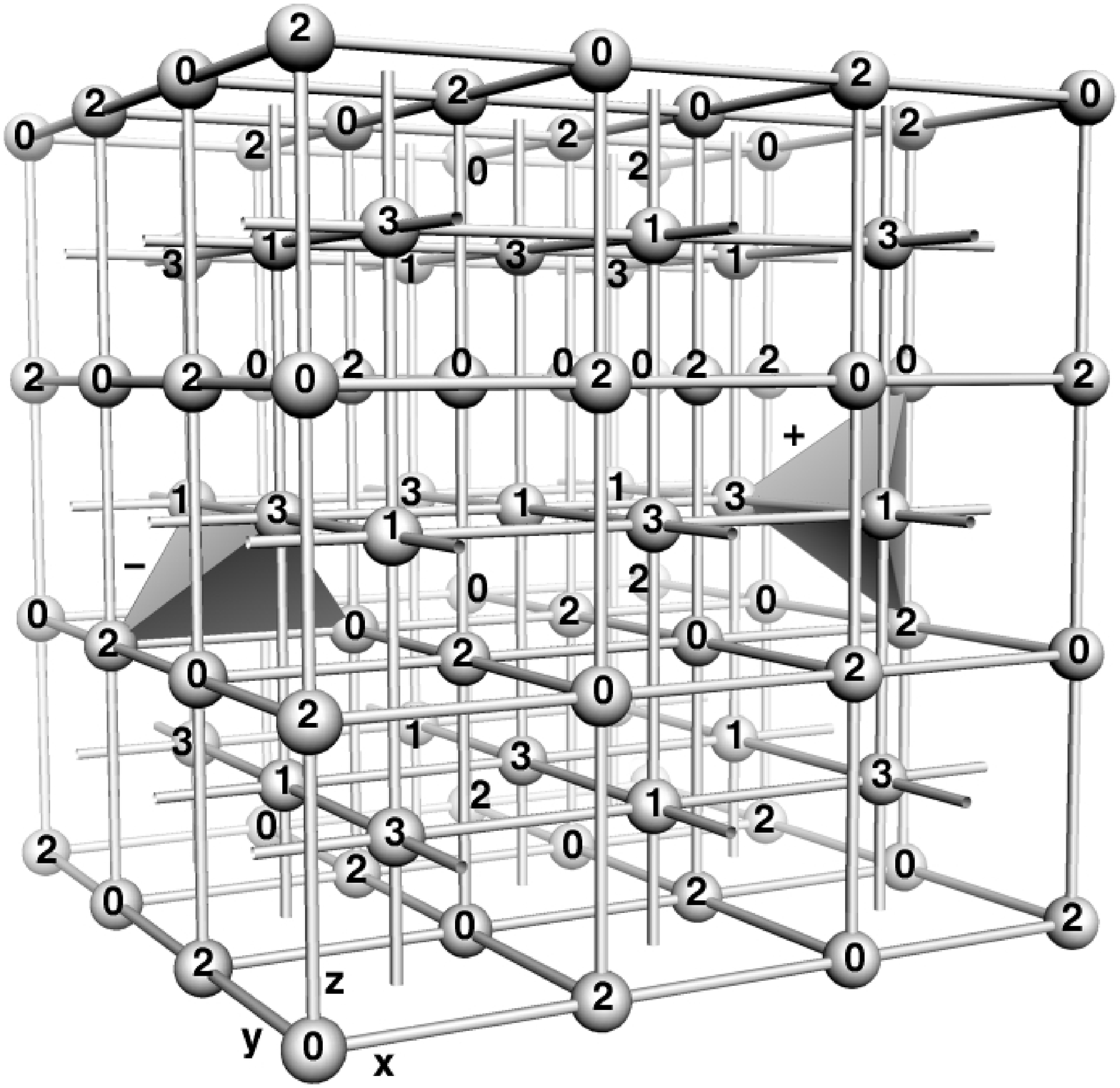,height=350pt}
\caption{%
The body-centered cubic lattice $D_3^\ast$.
The spheres colored $0$ to $3$ correspond to the four cosets
of the sublattice $D_3$.
The Delaunay cells fall into two classes, according as they can
be moved to coincide with the tetrahedron on the left or right, by a motion
that preserves the labeling of the spheres.
}\label{lattice-figure}
\end{figure}

These are the groups between $T_1$ and $N(T_1)$.
$N(T_1)$ is the automorphism group of the body centered cubic (bcc) lattice
indicated by the spheres
of Figure~\ref{lattice-figure} and $T_1$ is its odd subgroup.

The spheres of any one color $0$-$3$ correspond to a copy of the face-centered
cubic lattice
$$
\begin{array}{c}
	D_3=\{(x,y,z)\mid x,y,z \in {\mathbb Z},x+y+z \mbox{~is even}\}
\end{array}
$$
and those of the other colors to its cosets in the dual body-centered
cubic lattice
$$
\begin{array}{c}
	D_3^\ast=\{(x,y,z)\mid x,y,z\in {\mathbb Z} \mbox{~or~} x-\h,y-\h,z-\h
	\in {\mathbb Z}\}
\end{array}
$$
(the spheres of all colors).
We call these cosets $[0]$, $[1]$, $[2]$ and $[3]$, $[k]$ consisting
of the points for which $x+y+z\equiv \frac{k}{2} (\mbox{mod} 2)$.

Cells of the Delaunay complex are tetrahedra with one vertex of every color.
The tetrahedron whose center of gravity is $(x,y,z)$ is called {\em positive}
or {\em negative} according as $x+y+z$ is an integer $+\q$ or an integer
$-\q$. Tetrahedra of the two signs are mirror images of each other.

Any symmetry of the lattice
permutes the $4$ cosets and so yields a permutation $\pi$
of the four numbers $\{0,1,2,3\}$. We shall say that its image is
$+\pi$ or $-\pi$ according as it fixes or changes the signs of the
Delaunay tetrahedra.
The group $G_{16}$ of signed permutations so induced by all symmetries
of the lattice is $N(T_1)/T_1$.

\begin{figure}[c]
{
\psfragscanon
\psfrag{identity}{identity}
\psfrag{all symmetries}{all symmetries}
\psfrag{<(13)>}{$\langle(13)\rangle$}
\psfrag{<(01)(23)>}{$\langle(01)(23)\rangle$}
\psfrag{<(02)(13)>}{$\langle(02)(13)\rangle$}
\psfrag{<(13)_(02)>}{$\langle(13),(02)\rangle$}
\psfrag{<(01)(23)_(03)(12)>}{$\langle(01)(23),(03)(12)\rangle$}
\psfrag{<(0123)>}{$\langle(0123)\rangle$}
\begin{center}
\epsfig{file=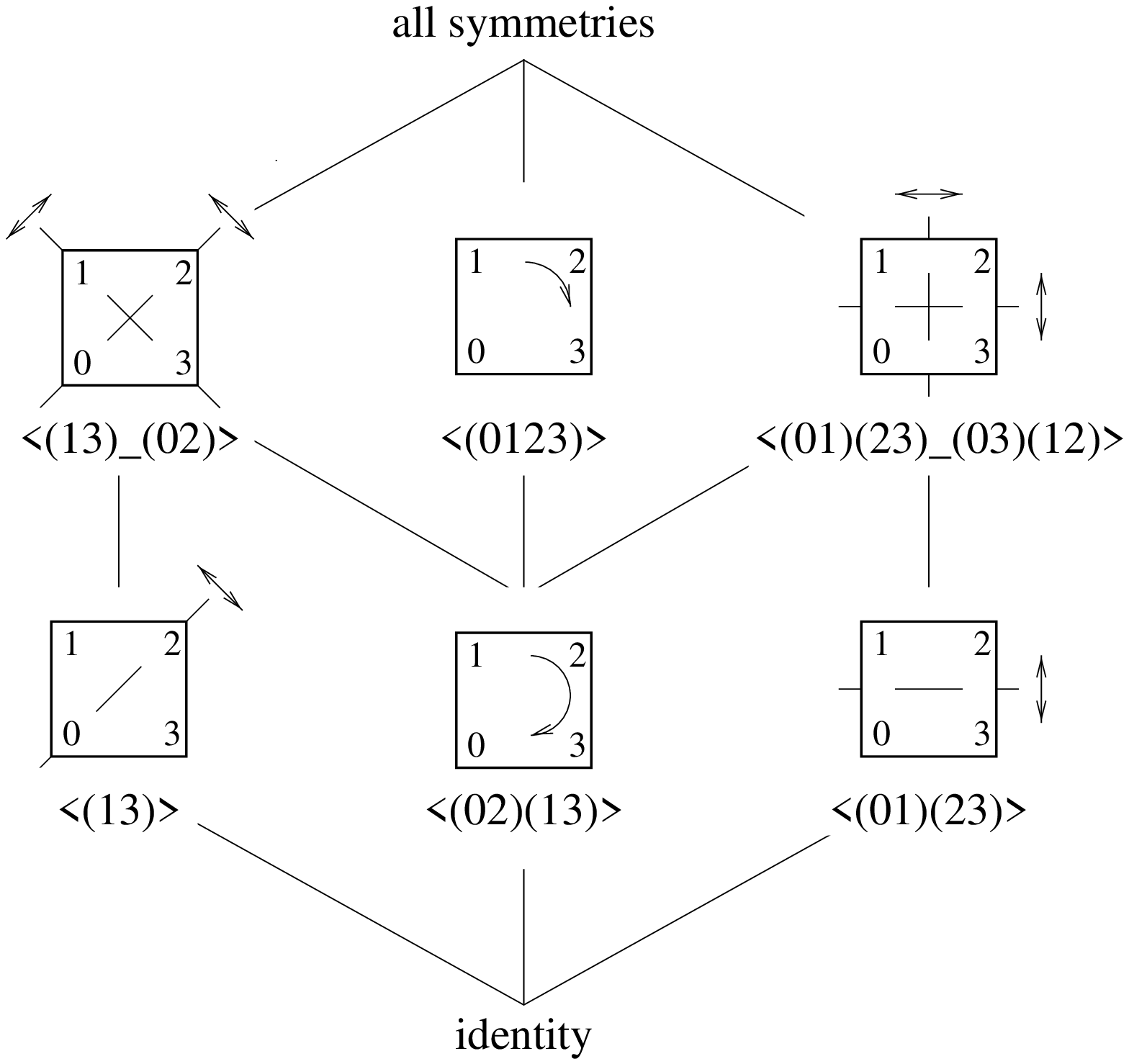,height=10cm}
\end{center}
\psfragscanoff
}
\caption{Subgroups of a dihedral group of order $8$.
The groups of order $2$ and $4$ on the left are generated by $1$ or
$2$ diagonal reflections; those on the right by $1$ or $2$ horizontal or
vertical reflections, and those in the center by a rotation of order $2$ or
$4$.}
\label{fig:subgroups}
\end{figure}

We obtain the 27 ``full'' space groups by selecting just those symmetries
that yield elements of some subgroup of $G_{16}$,
which is the direct product of the group $\{\PM 1\}$
of order $2$, with the dihedral group of order $8$ generated by the positive
permutations.
Figure~\ref{fig:subgroups} shows the subgroups of an abstract dihedral group
of order $8$, up to conjugacy.
Our name for a subgroup of order $n$ is
$n^{-}$, $n^{\o}$ or $n^{+}$, the superscript for a proper subgroup being
$\o$ for subgroups of the cyclic group of order $4$, and otherwise $-$
or $+$ according as elements outside of this are odd or even permutations.
The positive permutations form a dihedral group of order $8$, which
has 8 subgroups up to conjugacy, from which we obtain the following
$8$ space groups:
$$
\begin{array}{ccl}
8^{\o} &\mbox{from} &\{1,(02)(13),(0123),(3210),(13),(02),(01)(23),(03)(12)\}\\
4^{-} &\mbox{from} & \{1,(02)(13),(13),(02)\}\\
4^{\o} &\mbox{from} & \{1,(02)(13),(0123),(0321)\}\\
4^{+} &\mbox{from} & \{1,(02)(13),(01)(23),(03)(12)\}\\
2^{-} &\mbox{from} & \{1,(13)\} \mbox{~or~} \{1,(02)\}\\
2^{\o} &\mbox{from} & \{1,(02)(13)\}\\
2^{+} &\mbox{from} & \{1,(01)(23)\} \mbox{~or~} \{1,(03)(12)\}\\
1^{\o} &\mbox{from} & \{1\}.
\end{array}
$$

The elements of a
subgroup of $G_{16}$ that contains $-1$ come in pairs
$\PM g$, where $g$ ranges over one of the above groups, so we obtain
$8$ more space groups:
$$
\begin{array}{ccl}
8^{\o}:2 &\mbox{from}&\{\PM1,\PM(02)(13),\PM(0123),\PM(3210),\PM(13),\PM(02),\PM(01)(23),\PM(03)(12)\}\\
4^{-}:2 &\mbox{from}& \{\PM1,\PM(02)(13),\PM(13),\PM(02)\}\\
4^{\o}:2 &\mbox{from}& \{\PM1,\PM (02)(13),\PM(0123),\PM (0321)\}\\
4^{+}:2 &\mbox{from}& \{\PM1,\PM(02)(13),\PM(01)(23),\PM(03)(12)\}\\
2^{-}:2 &\mbox{from}& \{\PM1,\PM(13)\} \mbox{~or~} \{\PM1,\PM(02)\}\\
2^{\o}:2 &\mbox{from}& \{\PM1,\PM(02)(13)\}\\
2^{+}:2 &\mbox{from}& \{\PM1,\PM(01)(23)\} \mbox{~or~} \{\PM1,\PM(03)(12)\}\\
1^{\o}:2 &\mbox{from}& \{\PM1\}.
\end{array}
$$
Each remaining subgroup of $G_{16}$ is obtained by 
affixing signs to the elements of a certain
group $G$, the sign being $\p$ just for elements in some index 2 subgroup $H$
of $G$. If $H=N^i$, $G=2N^j$, we use the notation
$2N^{ij}$ for this. In this way we obtain 11 further space groups:
$$
\begin{array}{ccl}
8^{-\o} &\mbox{from}&\{
+1,+(02)(13),+(13),+(02),-(0123),-(3210),-(01)(23),-(03)(12)\}\\
8^{\o\o} &\mbox{from}&\{
+1,+(02)(13),+(0123),+(3210),-(13),-(02),-(01)(23),-(03)(12)\}\\
8^{+\o} &\mbox{from}&\{
+1,+(02)(13),+(01)(23),+(03)(12),-(13),-(02),-(0123),-(3210)\}\\
4^{--} &\mbox{from}&\{ +1, +(13), -(02)(13), -(02)\}\\
4^{\o-} &\mbox{from}&\{ +1, +(02)(13), -(13), -(02)\}\\
4^{\o\o} &\mbox{from}&\{ +1,+(02)(13), -(0123), -(0321)\}\\
4^{\o+} &\mbox{from}&\{ +1, +(02)(13), -(01)(23), -(03)(12)\}\\
4^{++} &\mbox{from}&\{ +1, +(01)(23), -(02)(13), -(03)(12)\}\\
2^{\o-} &\mbox{from}&\{ +1, -(13)\}\\
2^{\o\o} &\mbox{from}&\{ +1, -(02)(13)\}\\ 
2^{\o+} &\mbox{from}&\{ +1, -(01)(23)\}.\\
\end{array}
$$
(For these we have only indicated one representative of the conjugacy class.)

\medskip
Summary: A typical ``full group'' consists of all the symmetries
of Figure~\ref{lattice-figure} that induce a given group of signed
permutations.

\subsection{The 8 ``Quarter Groups''}

\begin{figure}[c]
\hfil
\epsfig{file=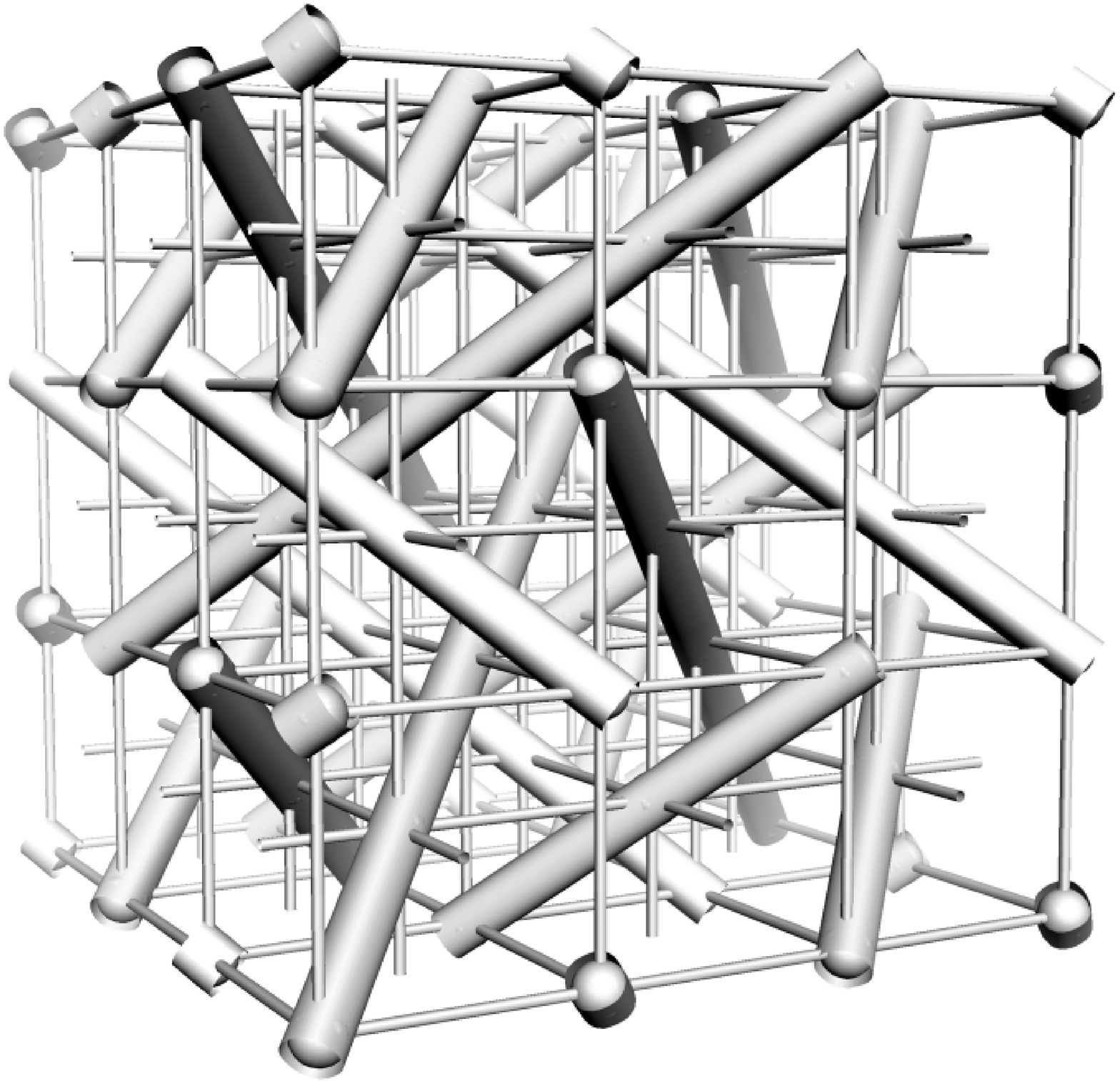,height=350pt}
\hfil
\caption{The cylinders represent the set of diagonal lines fixed by the eight
quarter groups.}
\label{lattice-line-figure}
\end{figure}

We obtain Figure~\ref{lattice-line-figure} from Figure~\ref{lattice-figure}
by inserting certain diagonal lines joining the spheres.

[The rules are that the sphere $(x,y,z)$
$$
\begin{array}{ccccc}
\multicolumn{4}{c}{\mbox{lies on a line in direction:}}\\
(1,1,1) & (-1,1,1) & (1,-1,1) & (1,1,-1)\\
\multicolumn{4}{c}{\mbox{according as:}}\\
x\equiv y\equiv z & y\not\equiv z\equiv x & z\not\equiv x\equiv y &
x\not\equiv y\equiv z & (\mbox{mod} 2)\\
\multicolumn{4}{c}{\mbox{if $x,y,z$ are integers, or according as:}}\\
x\equiv y\equiv z & z\not\equiv x\equiv y & x\not\equiv y\equiv z
& y\not\equiv z\equiv x
& (\mbox{mod} 2)\\
\multicolumn{4}{c}{\mbox{if they are not.]}}\\
\end{array}
$$

The automorphisms that preserve this set of diagonal lines form
the group $N(T_2)$ whose odd subgroup is $T_2$.
It turns out that these
automorphisms preserve each of the two classes of tetrahedra
in Figure~\ref{lattice-figure}
so that $N(T_2)/T_2$ is the order $8$ dihedral group
$G_8$ of positive permutations.

So the eight ``quarter groups'' $N^i/4$ (of index $4$ in the corresponding
group $N^i$) are obtained from the subgroups
of $G_8$ according to the following scheme:
$$
\begin{array}{ccl}
8^{\o}/4 &\mbox{from} &\{1,(02)(13),(0123),(3210),(13),(02),(01)(23),(03)(12)\}\\
4^{-}/4 &\mbox{from} & \{1,(02)(13),(13),(02)\}\\
4^{\o}/4 &\mbox{from} & \{1,(02)(13),(0123),(0321)\}\\
4^{+}/4 &\mbox{from} & \{1,(02)(13),(01)(23),(03)(12)\}\\
2^{-}/4 &\mbox{from} & \{1,(13)\} \mbox{~or~} \{1,(02)\}\\
2^{\o}/4 &\mbox{from} & \{1,(02)(13)\}\\
2^{+}/4 &\mbox{from} & \{1,(01)(23)\} \mbox{~or~} \{1,(03)(12)\}\\
1^{\o}/4 &\mbox{from} & \{1\}.
\end{array}
\old{
\begin{array}{ccl}
8^{\o}/4 &\mbox{from} &\{1,(0123),(02)(13),(3210),(13),(02),(01)(23),(03)(12)\}\\
4^{-}/4 &\mbox{from} & \{1,(13),(02),(02)(13)\}\\
4^{\o}/4 &\mbox{from} & \{1,(0123), (02)(13), (0321)\}\\
4^{+}/4 &\mbox{from} & \{1,(01)(23),(03)(12),(02)(13)\}\\
2^{-}/4 &\mbox{from} & \{1,(13)\} \mbox{~or~} \{1,(02)\}\\
2^{\o}/4 &\mbox{from} & \{1,(02)(13)\}\\
2^{+}/4 &\mbox{from} & \{1,(01)(23)\} \mbox{~or~} \{1,(03)(12)\}\\
1^{\o}/4 &\mbox{from} & \{1\}.
\end{array}
}
$$

Summary: A typical ``quarter group'' consists of all
the symmetries of Figure~\ref{lattice-line-figure} that
induce a given group of permutations.

\subsection{Inclusions Between the 35 Irreducible Groups}

\begin{figure}[c]
{
\renewcommand{\p}{{\small +}}
\renewcommand{\m}{{\small -}}
\renewcommand{\o}{\circ}
\newcommand{\N}{\hfill}
\newcommand{\Z}[1]{\begin{tabular}{p{1.6cm}}#1\end{tabular}}
\newcommand{\bs}{$\hspace{-3pt}$\mbox{$/$}}

\psfragscanon
\begin{center}
\epsfig{file=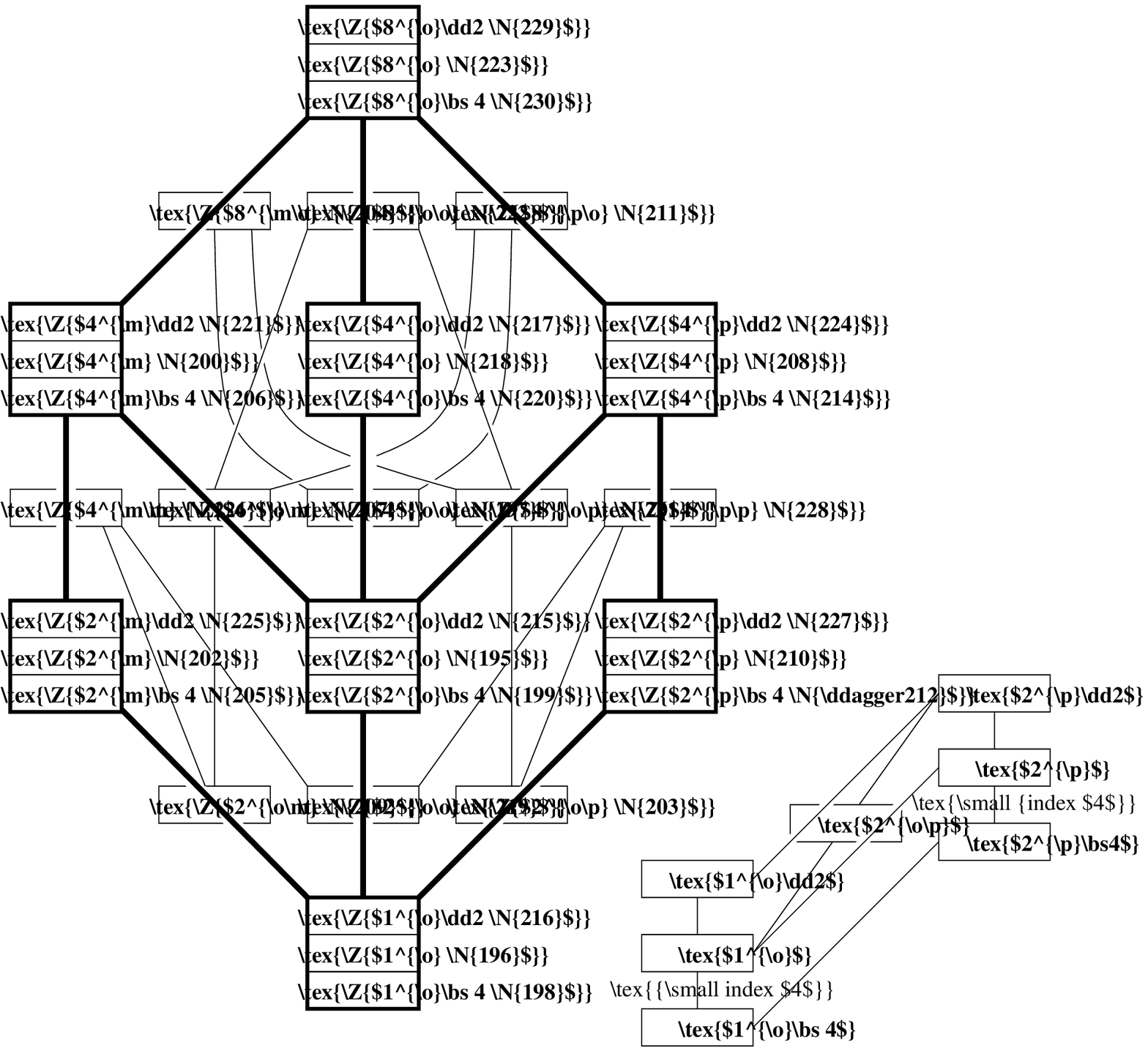,height=16.5cm}
\end{center}
\psfragscanoff
}
\caption{The 35 irreducible groups.
Each heavy edge represents several inclusions as in the inset.
With these conventions, the
figure illustrates all 83 minimal inclusions between these
groups: 2 for each of the 8 heavy boxes, 5 for each of the 11
heavy edges and 1 for each of the 12 thin edges.
(($\ddagger$): The group $2^+/4$ has two enantiomorphous forms,
with IT numbers 212 and 213.)
}\label{irreducible-figure}
\end{figure}

Our notation makes most of the inclusions obvious:
in addition to the containments $N^i$ in $2N^{ij}$ in $2N^j\dd2$,
each $G/4$
is index $4$ in $G$, which is index $2$ in $G\dd2$, and these three groups
are index $2$ in another such triple just if the same holds for the
corresponding subgroups of the dihedral group of order $8$.
All other minimal containments have the form $N^{ij}$ in $2N^{kl}$
and are explicitly shown in Figure~\ref{irreducible-figure}.

\subsection{Correspondence with the Fibered Groups}

Each irreducible group contains a fibered group to index $3$, and
when we started this work we hoped to obtain a nice notation for the
irreducible groups using this idea. Eventually we decided that the odd
subgroup method was more illuminating. However, we indicate
these relations in Table~2a.
We give a group $G$ secondary names such as
$$
	[\mbox{Jack}]\dd3 \quad \mbox{~or~} \quad (\mbox{Jill})\dd6, \quad
	\mbox{~say,}
$$
to mean
that $[\mbox{Jack}]$ is the group obtained by fixing the $z$
direction (up to sign) and $(\mbox{Jill})$ that obtained 
by fixing all three directions (up to sign).

\section{The 184 Reducible Space Groups}\label{outline-section}

We now consider the space groups that preserve some direction up to sign,
and accordingly can be given an invariant fibration.

\subsection{On Orbifolds and Fibration}

The rest of the enumeration is based on the concept of fibered orbifolds.
The {\em orbifold} of such a group is ``the space divided by the group'':
that is to say, the quotient topological space whose points are the orbits
under the group \cite{Thurston80,Scott83}.

For our purposes, a {\em fibration} is a division of the space into a system of
parallel lines. A {\em fibered space group} is a space group together with
a fibration that is invariant under that group.
On division by the group, the fibration of the space becomes a fibration of
the orbifold, each fiber becoming either a circle or an interval.
We call this a {\em fibered orbifold}.

The concept of a fibered space group differs slightly from that of a
reducible space group.
The latter are those for which there exists at least one
invariant direction and to obtain a fibered space group is to
make a fixed choice of such an invariant direction.
This distinction leads to what we call the ``alias problem'' discussed in
Section~\ref{alias-section}.

Although this paper was inspired by the orbifold concept,
we did not need to consider the 219 orbifolds of space groups individually.
We hope to discuss their topology in a later paper.

\subsection{Fibered Space Groups and Euclidean Plane Groups}

Look along the invariant direction of a fibered space group and you will
see one of the 17 Euclidean plane groups!

We explain this in more detail and introduce some notation.
Taking the invariant direction to be $z$, the action of any element of
the space group $G$ has the form:
$$
g: (x,y,z) \mapsto \big (a(x,y),b(x,y),c\plusminus z\big)
$$
for some functions $a(x,y), b(x,y)$, some constant $c$, and some sign
$\pm$.

Ignoring the $z$-component gives us the action
$$
g_H: (x,y)\mapsto \big(a(x,y),b(x,y)\big)
$$
of the corresponding element of the plane group.

We will call $g_H$ the {\em horizontal part} of $g$ and say
that it is {\em coupled} with the {\em vertical part}
$$
g_V: (z)\mapsto (c\plusminus z).
$$

\subsection{Describing the Coupling}

The fibered space group $G$ is completely specified by describing
the coupling between the ``horizontal'' operations $g_H$ and the
``vertical'' ones, $g_V$,
for which we use the notations $c\p$ and $c\m$,
where $c\p$ and $c\m$ are the maps $z\mapsto c+z$ and $z\mapsto c-z$.
We say that $h_H$ is plus-coupled or minus-coupled according as it
couples to an element  $c\p$  or  $c\m$.

Geometrically, $c\p$ is a translation through distance $c$, while $c\m$
is the reflection in the horizontal plane at height $\h c$.
It is often useful to note that
by raising the origin through a distance $\frac{d}{2}$,
we can augment by a fixed amount, or ``reset'', the constants in all
$c\m$ operations while fixing all $c\p$ ones.

In fact, any given horizontal operation is coupled with infinitely
many different vertical operations, since the identity is.
We study this by letting $K$
denote the {\em kernel}, consisting of all the vertical operations that are
coupled with the identity horizontal operation $I$.
Let $k$ be the smallest positive number for which $k\p \in K$.
Then the elements $nk\p$ $(n\in{\mathbb Z})$ are also in $K$.
If $K$ consists precisely of these elements, then the generic
fiber is a circle and we speak of a {\em circular fibration} and
indicate this by using $(\,)$'s in our name for the group.
If there exists some $c\m\in K$ then we can
suppose that $0\m \in K$, and then $K$ consists precisely of all $nk\p$ and
$nk\m$ $(n\in {\mathbb Z})$. In this case, the generic fiber is a closed
interval and we have an {\em interval fibration}, indicated
by using $[\,]$'s in the group name.

In our tables, we rescale to make $k=1$, so that $K$ consists
{\em either} of all elements $n\p$ (for a circular fibration) {\em or}
all elements $n\p$ and $n\m$ (for an interval one) for all integers
$n$.

\subsection{Enumerating the Fibrations}

To specify the typical fibration over a given plane group
$H=\<P,Q,R,\dots\>$, we merely have to assign vertical operations
$p\pm, q\pm ,r\pm \dots$, one for each of the generators $P,Q,R,\dots$.

The condition for an assignment to work is that it be a homomorphism,
modulo $K$. 
This means in particular that the vertical elements $c\pm$ need only be
specified modulo $K$, which allows us to suppose $0\leq c<1$,
which for a circular fibration is enough to make these elements unique.

For interval fibrations the situation is simpler, since we need only use
the vertical elements $0\p$ and $\h\p$.
(For since $0\m$ is in $K$ we can make the sign be $\p$; but also
$ c\p \cdot K \cdot (c\p)^{-1}=K $,
and since $(c\p) 0\m (c\p)^{-1}$ is the map taking $z$ to $2c-z$ this shows
that $2c$ must be an integer.)

\medskip
{\bf Example: the plane group
$H=632\cong \<\gamma,\delta,\epsilon\mid 1=\gamma^6=\delta^3=\epsilon^2=
\gamma\delta\epsilon\>$.}
\medskip

Here we take the corresponding vertical elements to be $c\pm, d\pm, e\pm$.
Then in view of $\gamma\delta\epsilon=1$,
it suffices to compute $d\pm$ and $e\pm$ .
The condition $\delta^3=1$ entails that the element
$d\pm$ is either $0\p$, $\t \p$, or $\T \p$.
(If the sign were $\m$, then the order of this element would have to be even.)
The condition $\epsilon^2=1$ shows that $e\pm$ can only be one of
$0\p$, $\h\p$, or $0\m$ and, (since the general case $e\m$
can be reset to
$0\m$), we get at most 9 circular fibrations, which
reduce to six by symmetries, namely $\gamma,\delta,\epsilon$ couple to
one of:
$$
\begin{array}{ccccc}
0\p 0\p 0\p &\quad & \T \p \t \p 0\p
&\cong& \t \p \T \p 0\p\\
\h\p 0\p \h\p &\quad & \frac{1}{6}\p \t \p \h\p
&\cong& \frac{5}{6}\p \T \p \h\p\\
0\m 0\p 0\m &\quad & \t \m \t \p 0\m
&\cong& \T \m \T \p 0\m.\\
\end{array}
$$

For interval fibrations, where we can only use
$0\p$ and $\h\p$, $\delta$ can only couple with $0\p$ (because
$\delta^3=1$), so we get at most $2$  possibilities,
namely that $\gamma,\delta,\epsilon,I$ couple with:
$$
\begin{array}{ccc}
	0\p0\p0\p0\m &\mbox{~or~} &\h\p0\p\h\p0\m.
\end{array}
$$

Studying these fibrations involves calculations
with products of the maps $k\pm$. Our notation makes this easy:
for example, the product $(a\m)(b\p)(c\m)$ is the map that takes
$z$ to $a-(b+(c-z))=(a-b-c)+z$; so $(a\m)(b\p)(c\m)=(a-b-c)\p$;
also, the inverse of $c\p$ is $(-c)\p$, while $c\m$ is its own inverse.

For example for the circular fibrations of $632$
we had $d\pm=0\p,\t \p, \T \p$ and
$e\pm= 0\p, \h\p, 0\p$ and these define $c\pm$ via the
relation $c\pm d\pm e\pm=0\p$.
So, if $d\pm=\t \p$ and $e\pm=\h\p$, then
$c\pm$ must be $c\p$, and since
$$
\begin{array}{c}
	c\p \t \p \h\p = (c+\t +\h)\p,
\end{array}
$$
$c$ must be $\frac{-5}{6}$, which we can replace by $\frac{1}{6}$.

The indicated isomorphisms are consequences of the isomorphism
that changes the sign of $z$, which allows us to replace $c,d,e$ by their
negatives (modulo $1$).
So we see that we have at most $6+2=8$ fibrations over the plane
group $632$. The ideas of the following section show that they are all
distinct.

\section{The Fibrifold Notation (for simple embellishments)}
\label{fibrifold-section}

Plainly, we need some kind of invariant to tell us when fibrations
really are distinct. A notation that corresponds exactly to the maps
$k\pm$ will not be adequate because they are far from being invariant:
for example, $0\m0\p0\m$ is equivalent to  $k\m 0\p k\m$ for every $k$.
We shall use what we call the {\em fibrifold notation},
because it is an invariant of the {\em fibered orbifold} rather than
the orbifold itself.
It is an extension of the orbifold notation that solved this problem
in the two-dimensional case \cite{Conway90,ConwayHuson99}.

The exact values of these maps are not, and should not be,
specified by the notation.
Rather, everything in the notation is an invariant of them up to
continuous variation (isotopy) of the group.
We usually prove this by showing how it corresponds to some feature of
the fibered orbifold.

We obtain the fibrifold notation for a space group $G$ by
``embellishing'', or adding information
to the orbifold notation for the two-dimensional group $H$ that is its
horizontal part. (Often the embellishment consists of doing nothing.
In this section we handle only the simple embellishments that
can be detected by local inspection of the fibration.)

\subsection{Embellishing a ring symbol $\o$}

A ring symbol corresponds to the relations
$\up{\alpha}\o\up{X}\up{Y}$: $\alpha=[X,Y]$.
We embellish it to $\o$ or $\bar{\o}$ to mean that $X,Y$ are both plus-coupled
or both minus-coupled respectively.
[We can suppose that $X$ and $Y$ couple to the same sign in view of the
three-fold symmetry revealed by adding a new generator
$Z$ satisfying the relations $XYZ=1$ and $X^{-1}Y^{-1}Z^{-1}=\alpha$.]

This embellishment is a feature of the fibration since it tells
us whether or not the fibers have a consistent orientation over
the corresponding handle.

\subsection{Embellishing a gyration symbol $G$}

Such a symbol corresponds to relations
$\up{\gamma}G$: $\gamma^G=1$. We embellish it to
$G$ or $G_g$ according as $\gamma$ is minus-coupled or coupled
to $\frac{g}{G}\p$.
The relation $\gamma^G=1$ implies that $G$ must be even if $\gamma$ is
minus-coupled and that $c$ must be a multiple of $\frac{1}{G}$, if
$\gamma$ couples to $c\p$.

The embellishment is a feature of the fibration because
the behavior of the latter
at the corresponding cone point determines an action of
the cyclic group of order $G$ on the circle.

\subsection{Embellishing a kaleidescope symbol $\s A B \dots C$}

Here the relations are
$\up{\lambda}\s\up{P}A\up{Q}B\dots\up{R}C\up{S}$:
$$
    1=P^2=(PQ)^A=Q^2=\dots= R^2=(RS)^C=S^2 \mbox{~and~}
	\lambda^{-1}P\lambda =S.
$$
We embellish $\s$ to $\s$ or $\bar{\s}$
according as $\lambda$ is plus-coupled or minus-coupled.
This embellishment is an fibration feature because it tells us whether
the fibers have or do not have a consistent orientation when
we restrict to a deleted neighborhood of a string of mirrors.

The coupling of Latin generators is indicated by inserting
$0$, $1$ or $2$ dots into the corresponding spaces of the orbifold name.
Thus $AB$ will be embellished to $AB$, $A\d B$ or $A \dd B$ according
as $Q$ is minus-coupled or coupled to $0\p$ or $\h\p$.
These embellishments tell us how the fibration behaves above the
corresponding line segment. The generic fiber is identified with
itself by a homomorphism whose square is the identity, namely
$$
\begin{array}{cccll}
	\d & (\mbox{e.g.~} A\d B) &
		\mbox{~corresponds to a mirror~} & \phi(x)=x\\
	\dd & (\mbox{e.g.~} A \dd B) &
		\mbox{~is a M\"obius map~} & \phi(x)=x+\h,& \mbox{~and}\\
	\mbox{a blank} & (\mbox{e.g.~} A\, B)&
		\mbox{~corresponds to a link~} & \phi(x)=-x.\\
\end{array}
$$

The numbers $A,B,\dots,C$ are the orders of the products like $PQ$
for adjacent generators. So, as in the gyration case, we
embellish $A$ to $A$ or $A_a$ according as $PQ$ is minus-coupled
or coupled to $\frac{a}{A}\p$.
There are two cases:
If $P$ and $Q$ are both minus-coupled, say to $p\m$ and $q\m$, then
$PQ$ is coupled to $(p-q)\p$, and so that $\frac{a}{A}$ is $p-q$.
If they are both plus-coupled,
then the subscript $a$ is determined as $0$ or $\frac{A}{2}$ by whether the
punctuation marks $\d$ or $\dd$ on each side of $A$
are the same or different, and so we often omit it.

In pictures we use corresponding embellishments
\begin{center}
{
\def\VLINE{\mbox{\LARGE $|$}}
\begin{tabular}{p{1cm}p{1cm}p{1cm}p{1cm}p{1cm}}
	\hfil $\VLINE\d$  \hfil & \hfil $\VLINE \dd$ \hfil & \hfil $\VLINE$ \hfil
	& \hfil $n_d$ \hfil & \hfil $n$ \hfil\\
\end{tabular}
}
\end{center}
to show that the appropriate reflections or rotations are coupled to:
\begin{center}
\begin{tabular}{p{1cm}p{1cm}p{1cm}p{1cm}p{1cm}}
	\hfil $0\p$ \hfil & \hfil $\h \p$ \hfil & \hfil $k\m$ \hfil &\hfil
	$\frac{d}{n}\p$ \hfil & \hfil $k\m$. \hfil \\
\end{tabular}
\end{center}

\subsection{Embellishing a cross symbol $\x$}

Here the relations are $\up{\omega}\x\up{Z}$: $Z^2=\omega$, and
we embellish $\x$ to $\x$ or $\bar{\x}$ according as $Z$ is plus-coupled
or minus-coupled (indicating whether the fibers do or do not have
a consistent orientation over the corresponding crosscap).

\subsection{Other embellishments exist!}

These simple embellishments should suffice for a first reading,
since they suffice for many groups.
In a Section~\ref{completing-section} we shall describe more subtle ones that
are sometimes
required to complete the notation.

\section{The Enumeration: Detecting Equivalences}

The enumeration proceeds by assigning coupling maps $k\pm$
to the operations in all possible ways that yield homomorphisms
modulo $K$, and describing these in the (completed) fibrifold notation.
What this notation captures is exactly the isotopy class of a given
fibration over a given plane crystallographic group $H$, in other words
the assignment of the maps $k\p$ and $k\m$ to all elements of $H$,
up to continuous variation.

\subsection{The Symmetry Problem}

\begin{figure}[c]
\hfil
\begin{tabular}{cc}
\epsfig{file=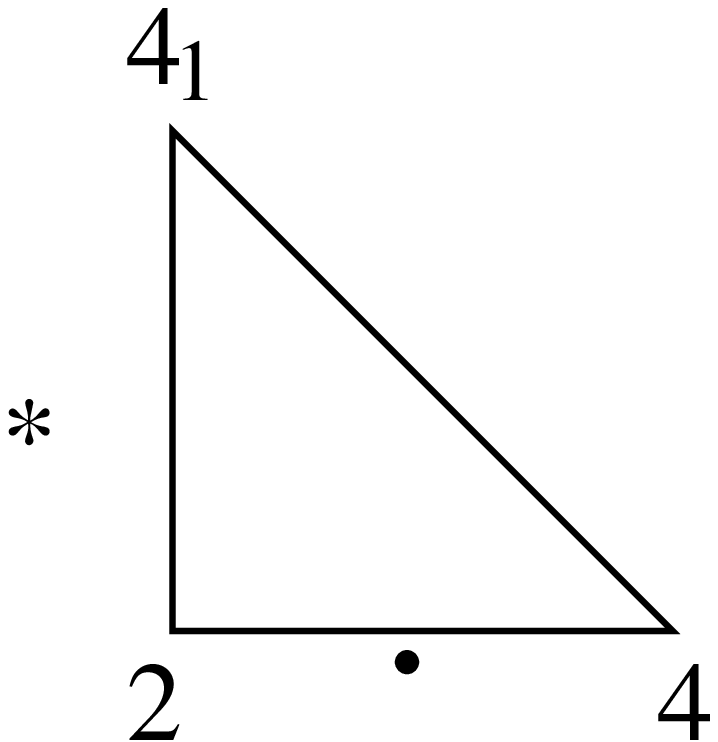,height=3cm} &
\epsfig{file=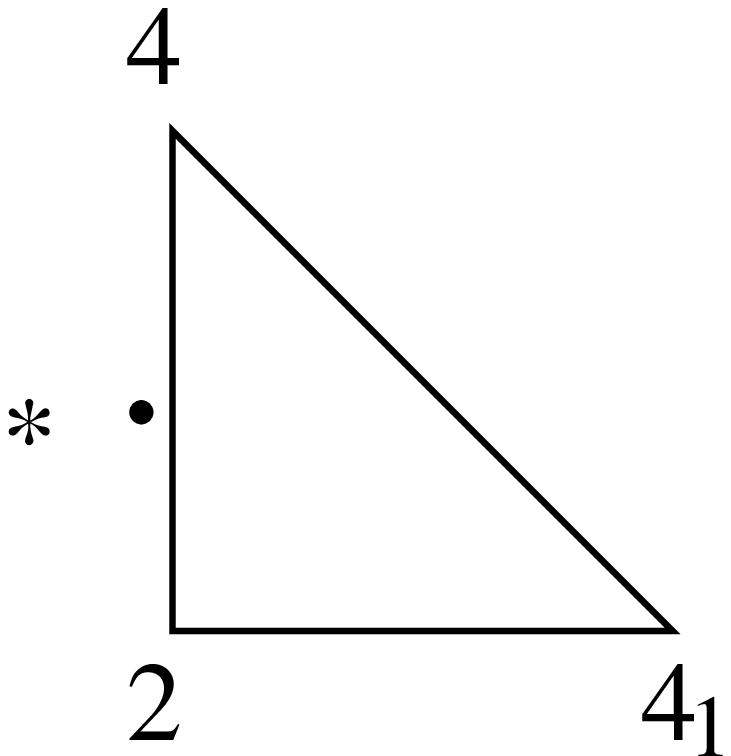,height=3cm}\\
$(\s 4_1 4 \d 2)$ & $(\s \d 4 4_1 2)$\\
\end{tabular}
\hfil
\caption{Two fibrations related by a symmetry.}
\label{fig:symmetry-example}
\end{figure}

The two fibrations of $\s\up{P}4\up{Q}4\up{R}2$ notated
$(\s 4_1 4 \d 2)$ and $(\s \d 4 4_1 2)$
are distinct in this sense - in the first $R$ maps to the identity,
while in the second $P$ does.
However, they still give a single three-dimensional group because the plane
group $H=\s442$ has a symmetry taking $P,Q,R$ to $R,Q,P$,
see Figure~\ref{fig:symmetry-example}

So our problem is really to enumerate isotopy classes of
fibrations up to symmetries: we discuss the different types of symmetry
in Appendix~III. Some of them are quite subtle -
how can we be sure that we have accounted for them all?

To be quite safe, we made use of the programs of Olaf Delgado and Daniel Huson
\cite{DelgadoHuson96} which compute various
invariants that distinguish three-dimensional space groups.
Delgado and Huson used these to find the 
``IT number'' that locates the group in the International
Crystallographic Tables \cite{Hahn83};
but our enumeration uses the invariants only, and so is logically
independent of the international tabulation.
The conclusion is that the reduced names we introduce in
Section~\ref{alias-section} account for all the equivalences
induced by symmetries between fibrations over the same plane group.

\subsection{The Alias Problem}\label{alias-section}

However, a three-dimensional space group $G$ may have several
invariant directions, which may correspond to different fibrations
over distinct plane groups $H$.
So we must ask: which sets of names - we call them {\em aliases} -
correspond to fibrations of the same group?
This can only happen when the programs of Delgado and Huson yield
the same IT number, a remark that does not in fact depend on the
international tabulation, since it just means that they are the only groups
with certain values of the invariants.

\begin{figure}[c]
\begin{center}
\begin{tabular}{cc}
\epsfig{file=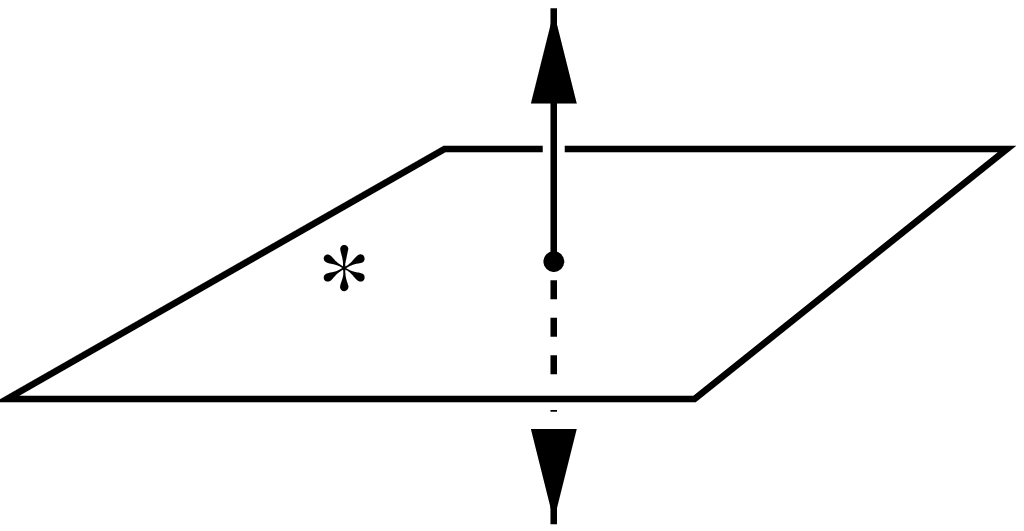,height=2cm} &
\epsfig{file=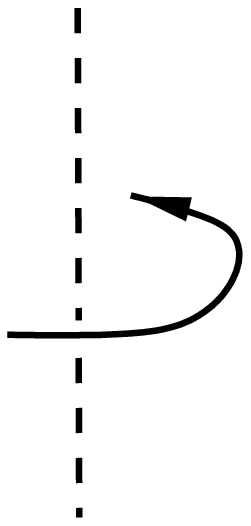,height=2cm} \\
(a) Reflection in $z=0$ & (b) half-turn about $z$-axis\\
\epsfig{file=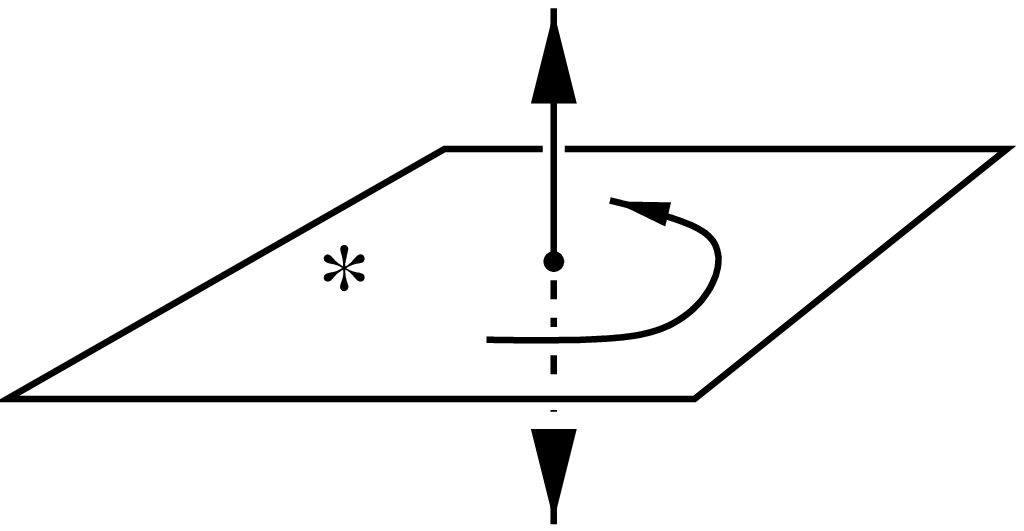,height=2cm} &
\epsfig{file=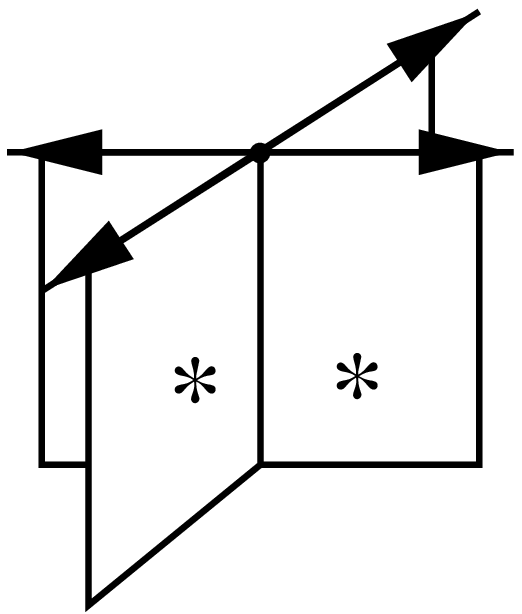,height=2cm} \\
(c) reflection and half-turn & (d) reflections in $x=0$ and $y=0$\\
\end{tabular}
\end{center}
\caption{The point groups $\s$, $22$, $2\s$ and $\s22$ each determine
a unique canonical direction.}
\label{directions-fig}
\end{figure}

The alias problem only arises for the reducible point groups, namely
$1$, $\x$, $\s$, $22$, $2\s$, $\s22$, $222$ and $\s222$.
In fact there is no problem for $1$ and $\x$, since for these cases
there is just one fibration. For $\s$, $22$, $2\s$ and $\s22$
the point group determines a unique canonical direction (see
Figure~\ref{directions-fig})
and so a unique {\em primary} name. Any other, {\em secondary} name,
must be an alias for some primary name, solving the alias problem
in these cases.

\begin{figure}[c]
\begin{center}
\begin{tabular}{cp{1cm}c}
\psfig{file=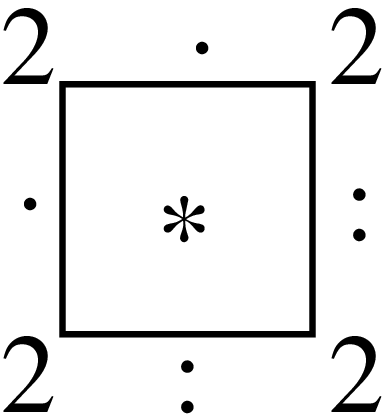,height=2cm} &  &
\psfig{file=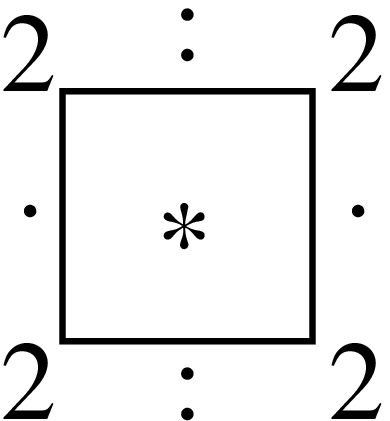,height=2cm}\\
(a) &  &(b)\\
\end{tabular}
\end{center}
\caption{Diagram (a) possesses a symmetry that interchanges the $x$ and
$y$ axes, whereas diagram (b) does not.
}\label{symmetric-fig}
\end{figure}

We are left with the point groups $\s222$ and $222$,
which always have three fibrations in orthogonal directions.
A permutation of the three axes might lead to an isotopic group.
If the number of such permutations is:
\begin{itemize}
\item[] 1, we have three distinct asymmetric names, say $\{A,A',A''\}$,
\item[] 2, we have two  distinct names, one symmetric and one asymmetric,
say $\{S,A\}$,
\item[] 3, we have just one asymmetric name, say $\{A\}$, or
\item[] 6, we have just one symmetric name, $\{S\}$,
\end{itemize}
where an axis, or the corresponding group name, is called {\em symmetric}
if there is a symmetry interchanging the other two axes.
This can be detected from the fibrifold notation,
for example, it is apparent from Figure~\ref{symmetric-fig}
that $\s\d2\d2\dd2\dd2$ is symmetrical and
$\s\d 2\dd 2\d 2\dd 2$ is not.

Table~2b provisionally lists all names that have the same invariant values
in a single line. We now show that the corresponding alias-sets
$$
	\{A,A',A''\}, \{S,A\}, \{A\}, \{S\}
$$
are correct. For if not, some $\{A,A',A''\}$ or $\{S,A\}$
would correspond to two or more groups, one of which would have a single
asymmetric name. But we show that there is only one such group,
whose name $(2_1 2\bar{\s}:)$ was not in fact in a set of type $\{A,A',A''\}$
or $\{S,A\}$.

A group $G$ with a single asymmetric name must be isotopic to
that obtained by cyclically permuting $x,y,z$: this isotopy
will become an automorphism if we suitably rescale the axes.
Adjoining this automorphism leads to an irreducible group in which $G$ has
index $3$. But inspection of Table~2a
reveals that the only asymmetric name
for which this happens is $(2_1 2\bar {\s}:)$.

So indeed the primary and secondary names in any given line
of Table~2b are aliases for the same group.
Our rules for selecting the primary name are:
\begin{enumerate}
\item
A unique name is the primary one.
\item
Of two names, the primary name is the symmetrical one.
\item
Otherwise, the primary name is that of a fibration over $22\s$.
\item
Finally, we prefer  $(2_0 2\bar{\s}\d)$ to $[2_1 2_1 \s\dd]$
and $(2_0 2\bar{\s}\dd)$ to $(2_1 2 \bar{\s}_1)$.
\end{enumerate}

These conventions work well because the three-name cases all involve a
fibration over $22\s$ (which was likely because the $x$ and $y$ axes
can be distinguished for $22\s$). The groups in the last rule are those
with two such fibrations.

\section{Completing the Embellishments}\label{completing-section}

We now ask what further embellishments are needed to specify a fibration
up to isotopy? We will find that it suffices to add subscripts $0$ or $1$
to some of the symbols $\o$, $\s$ and $\x$.
This section can be omitted on a first reading, since
for many groups in the tables these subtle embellishments are not needed.
We shall show as we introduce them that 
they determine the coupling maps up to isotopy,
so that no further embellishments are needed.

\subsection{Ring symbol $\o$}

Here if $X$ and $Y$ couple to $x\p$ and $y\p$, then the relation
$[X,Y]=\alpha$ shows that $\alpha$ is automatically coupled to $0\p$.
So the space groups obtained for arbitrary values of $x$ and $y$ are
isotopic and we need no further embellishment.

If, however,  there are one or more
embellishments to $\bar {\o}$, say those in the relations
$$
	\up{\alpha_1}\o\up{X}\up{Y}\dots \up{\alpha_n}\o\up{U}\up{V}
$$
then the coupling of the global relation will have the form
$$
\begin{array}{ccccc}
\underbrace{\alpha_1} & \dots &\underbrace{\alpha_n} &
\underbrace{(\gamma\dots\omega)} &=1\\
\downarrow & &\downarrow &\downarrow\\
2(x-y)\p & \dots & 2(u-v)\p & k+\\
\end{array}
$$
which restricts the variables $x,y,\dots,u,v$ only by the condition
that $2(x-y+\dots+u-v)$ be congruent to $-k$ (modulo 1).

So if $k$ is already determined this leaves just two values (modulo 1)
for $x-y+\dots+u-v$, namely $\frac{i-k}{2}$ ($i=0$ or $1$);
we distinguish when necessary by adding $i$ as a subscript. Once again,
the individual values of $x,y,\dots,u,v$ do not matter since they can be
continuously varied in any way that preserves the truth of
$$
	x-y+\dots+u-v=\frac{i-k}{2}.
$$

\subsection{Gyration symbol $G$}

A gyration $\gamma$ can only be coupled to $0\p$, $\h\p$ or $c\m$.
If we suppose that the minus-coupled gyrations are
$\gamma_1 \mapsto c_1\m, \dots, \gamma_n \mapsto c_n\m$, then
the values of the $c_i$ are unimportant for the local relations, while
the global relation involves only
$c_1-c_2+\dots\pm c_n$, and the $c_i$ can be varied in any way
that preserves this sum. So all solutions are isotopic and no
more subtle embellishment is needed.

\subsection{Kaleidescope symbol $\s A B \dots C$}

The simple embellishments suffice for the relations
$$
    1=P^2=(PQ)^A=Q^2=\dots=R^2=(RS)^C=S^2,
$$
so we need only discuss the relation $\lambda^{-1}P\lambda=S$
and the global relation.
We have already embellished $\s$ to $\s$ or $\bar{\s}$ according
as $\lambda$ couples to an element $l\p$ or $l\m$.

Specifying $l$ more precisely is difficult, because we need two rather
complicated conventions for circular fibrations, neither of which seems
appropriate for interval ones.

For interval fibrations we simply embellish $\s$ to $\s_0$ or $\s_1$
according as $\lambda$ couples to $0\p$ or $\h\p$ (the only two
possibilities).

What does the relation $\lambda^{-1}P\lambda=S$ tell us about the number $l$
in the circular-fibration case?

The answer turns out to be `nothing', if any one of the Latin generators
$P,\dots,S$ couples to a translation $k\p$.
To see this, it suffices to suppose that $P$ couples to
$0\p$ or $\h\p$,
and then $S$, being conjugate to $P$, must couple to the {\em same} translation
$0\p$ or $\h\p$.
But since these two translations are central they conjugate to themselves
by {\em any} element $l\pm$, so the relation 
$\lambda^{-1} P \lambda=S$ is automatically satisfied, and we need no
further embellishment.

The only remaining case is when all of $P,\dots,S$ couple to reflections
$p\m, q\m, \dots, s\m$.
In this case we have already embellished
$A\, B \dots C$ to $A_a B_b\dots C_c$ and we have
$$
	q=p-\frac{a}{A},\dots,s=r-\frac{c}{C},
$$
so
$$
	s=p-(\frac{a}{A}+\frac{b}{B}+\dots\frac{c}{C})
	=p- \Sigma, \mbox{~say.}
$$
Then according as $\lambda$ couples to $l\p$ or $l\m$, the relation
$\lambda^{-1} P \lambda=S$ now tells us that (modulo 1)
$$
-l+p-(l+z)=p-\Sigma-z \quad\mbox{~or~}\quad l-(p-(l-z))=p-\Sigma-z,
$$
whence (again modulo 1)
$$
	2l=\Sigma, \quad\mbox{~or~}\quad 2l=2p-\Sigma,
$$
and so finally
$$
l=\frac{\Sigma+i}{2}\quad\mbox{~or~}\quad p-l=\frac{\Sigma+i}{2} \quad
(i=0 \mbox{~or~} 1).
$$

In other words, our embellishments so far determine $2l$, but
not $l$ itself (modulo 1). We therefore further embellish $\s$ or
$\bar{\s}$ to $\s_i$ or $\bar{\s}_i$ to distinguish these two cases.
Their topological interpretation is rather complicated.
Two adjacent generators $P$ and $Q$ correspond to reflections of
the circle that each have two fixed points: let $p_0$ and $p_1$
be the fixed points for $P$.
Then if the product of the two reflections has rotation number $\frac{a}{A}$,
the fixed points for $Q$ are rotated by $\frac{a}{2A}$ and
$\frac{a+A}{2A}$ from $p_0$; call these $q_0$ and $q_1$ respectively.
In this way the fractions $\frac{a}{A},\frac{b}{B},\dots$
enable us to continue the naming of fixed points all around the circle.
The subscripts $0$ and $1$ tell us whether when we get back
to $p_0$ and $p_1$ they are in this order or the reverse.

\subsection{Cross symbol $\x$}

We have embellished a cross symbol to $\x$ or $\bar {\x}$ according
as $Z$ couples to $z\p$ or $z\m$.

What about the value of $z$?
The only relations involving $Z$ and $\omega$ are
$$
	Z^2=\omega \quad \mbox{and} \quad 1=\alpha\dots\omega.
$$
In the case $\bar{\x}$ the first of these implies that $\omega\mapsto 0\p$ for
any $z$, and so these relations will remain true if $z$ is continuously
varied. In other words, the space groups we obtain here for different
values of $z$ are mutually isotopic.

When there are $\x$ symbols, say
$\up{\psi}\x\up{Y}\cdots\up{\omega}\x\up{Z}$, not embellished to $\bar{\x}$ the
situation is different, because we have
$$
	Y\mapsto y\p,\dots, Z\mapsto z\p
$$
and so the relations $Y^2=\psi,\dots, Z^2=\omega$ imply
$$
	\psi \mapsto 2y\p,\dots,\omega\mapsto 2z\p
$$
and now the coupling of the global relation must take the form
$$
\begin{array}{ccccc}
\underbrace{(\gamma\dots\lambda\dots)} &  \underbrace{\psi} & \dots&
 \underbrace{\omega} &=1\\
\downarrow & \downarrow & & \downarrow\\
k\p & 2y\p & \dots & 2z\p\\
\end{array}
$$
showing that $2(y+ \dots +z)$ must be congruent to $-k$ (modulo 1), where
$k$ may be already determined.
A subscript $i=0$ or $1$ on such a string of $\x$ symbols
will indicate that $y+\dots+z=\frac{-k+i}{2}$.

However, it may be that the
remaining relations allow $k$ to be continuously varied, in which case
no subscript is necessary.

\section*{Appendix I: Crystallographic Groups and Bieberbach's Theorem}

\old{
A {\em crystallographic group} is a cocompact discrete group
of isometries of Euclidean space. The following
result is fundamental for the theory of crystallographic groups;
it was first proven by Schoenfliess in 1891 in three dimensions and then
by Bieberbach in 1911 in the $d$-dimensional case:
{\em
Any crystallographic group $G$ contains a subgroup $T$ of
finite index consisting of parallel translations.}
The quotient group $G/T$ is called the {\em point group} of $G$.

In the three-dimensional case,
the group of translations is ${\mathbb Z}^3$ and has finite index in
$G$. Further, the point group is a finite subgroup of
$\mbox{GL}(3,{\mathbb Z})$, and is the fundamental group of a
two-dimensional spherical group. It maps injectively to
$GL(3,{\mathbb Z}/n)$, for $n\geq 3$, from which one can readily determine the
possible point groups.
}

A {\em $d$-dimensional crystallographic group $G$} is a discrete co-compact
group of isometries of $d$-dimensional Euclidean space ${\mathbb E}^d$.
In other words, $G$ consists of isometries (distance preserving maps),
any compact region contains at most finitely many $G$-images of a
given point and the $G$-images of some compact region cover ${\mathbb E}^d$.

Every isometry of ${\mathbb E}^d$ can be written as a pair $(A,v)$ consisting
of a
$d$-dimensional orthogonal matrix  $A$ and a $d$-dimensional vector $v$.  The
result of applying $(A,v)$ to $w$ is $(A,v)w:=Av+w$ and thus the product of
two pairs $(A,v)$ and $(B,w)$ is $(AB,Aw+v)$. The matrix $A$ does not
depend on the choice of origin in ${\mathbb E}^d$.

An isometry $(A,v)$ is a pure translation exactly if $A$ is the identity
matrix.  The translations in a given isometry group $G$ form a normal
subgroup $T(G)$, namely the kernel of the homomorphism
$$
\begin{array}{rccl}
        \rho\colon& G &\to&     O(d)\\
        & (A,v)        &\mapsto& A,
\end{array}
$$
and so there is an exact sequence
$$
        0 \longrightarrow T(G)
          \stackrel{\subseteq}{\longrightarrow} G
          \stackrel{\rho}{\longrightarrow} \rho(G) \to 0.
$$

The image $\rho(G)$ is called the point group of $G$. Since $T(G)$ is
abelian, the conjugation action of $G$ on $T(G)$ is preserved by $\rho$, 
and the point group acts in a natural way on the group of translations.

If $T(G)$ spans (i.e.\ contains a base for) ${\mathbb R}^d$, then $T(G)$ is a maximal
abelian subgroup of $G$. To see this, consider an element $a=(I,v)$ in
$T(G)$ and an element $b=(A,w)$ in $G$. If $a$ and $b$ commute, we have
$w+v=Av+w$, so $v=Av$. Assume that $b$ commutes with every element in
$T(G)$. Since $T(G)$ spans ${\mathbb R}^d$, it follows that $A$ is the identity
matrix and that $b$ is a translation.

The translations subgroup of a crystallographic group is discrete and
therefore isomorphic to ${\mathbb Z}^m$ for some $m\leq d$. A famous theorem by
Bieberbach published in 1911 states (among other things)
that for a crystallographic group $G$, the point group $\rho(G)$ is finite
and $T(G)$ has full rank. Consequently, $\rho(G)$ is isomorphic to a finite
subgroup of $\mbox{GL}(d,{\mathbb Z})$.

Bieberbach also proved that in each dimension, there are only finitely many
crystallographic groups and that any two such groups are abstractly isomorphic
if and only if they are conjugate to each other by an affine map. For a
modern proof of Bieberbach's results, see \cite{Thurston97}.

The point group of a $3$-dimensional crystallographic group is a finite
subgroup of $O(3)$ and thus the fundamental group of a spherical
$2$-dimensional orbifold, which contains no rotation of order other than
$2$, $3$, $4$ or $6$ (the so-called {\em crystallographic restriction}).
All such groups are readily enumerated, for example using the two-dimensional
orbifold notation \cite{ConwayHuson99}.

\section*{Appendix II: On Elements of Order 3}

We show that any irreducible space group contains elements of order $3$.
For the five irreducible point groups all contain $332$, which is generated
by the two operations:
$$
	r: (x,y ,z) \mapsto (y,z,x) \mbox{~and~} s: (x,y,z) \mapsto (x,-y,-z).
$$
We may therefore suppose that our space group contains the operations:
$$
	R: (x,y,z) \mapsto (a+y,b+z,c+x) \mbox{~and~}
	S: (x,z,y) \mapsto (d+x,e-y,f-z).
$$

We now compute the product $SR^3S^{-1}R$:
$$
\begin{array}{c}
(x,y,z)\\
\, \downarrow S\\
(d+x,e-y,f-z)\\
\, \downarrow R^3\\
(d+x+a+b+c,e-y+a+b+c,f-z+a+b+c)\\
\, \downarrow S^{-1}\\
(x+a+b+c,y-a-b-c,z-a-b-c)\\
\, \downarrow R\\
(y-b-c,z-a-c,x+a+b+2c).
\end{array}
$$
Modulo translations, this becomes the order $3$ element
$(x,y,z)\mapsto (y,z,x)$ of the point group, but it has a fixed
point, namely $(-b-c,0,a+c)$, and so is itself of order $3$.

Now we show that there are only two possibilities for the odd subgroup $T$
of an irreducible space group. The axes of the order $3$ elements
are in the four directions
$$
	(1,1,1),\, (1,-1,-1),\, (-1,1,-1),\, \mbox{~and~} \,(-1,-1,1).
$$

The simplest case is when no two axes intersect.
In this case we select a closest pair of non-parallel axes, and
since any two such pairs are geometrically similar, we may suppose that
these are the lines in directions
$$
	(1,1,1) \mbox{~through~} (0,0,0) \mbox{~and~}
	(1,-1,-1) \mbox{~through~} (0,1,0),
$$
whose shortest distance is $\sqrt{\h}$.
Now the rotations of order three about these two lines generate
the group $T_2$, whose order $3$ axes - we call them the ``old'' axes -
are shown in Figure~\ref{lattice-line-figure}.
We now show that $T=T_2$; for if not, there must be another order $3$ axis
not intersecting any of the ones of $T_2$. However
(see Figure~\ref{lattice-line-figure}),
the entire space is partitioned into $\h\times \h \times \h$ cubes
which each have a pair of opposite vertices on axes of $T_2$,
and any such cube is covered by the two spheres of radius $\sqrt\h$
around these opposite vertices.
Any ``new'' axis must intersect one of these cubes, and therefore its distance
from one of the old axes must be less than $\sqrt{\h}$, a contradiction.

If two axes intersect, then each axis will contain infinitely many such
intersection points; we may suppose that a closest pair of these 
are the points $(0,0,0)$ and $(\h,\h,\h)$, in the direction $(1,1,1)$.
Then the group contains the order $3$ rotations about the four axes
through each of these points, which generate $T_1$, the odd subgroup
of the body centered cubic (bcc) lattice.

We call the order $3$ axes of $T_1$ the ``old'' axes: they consist of
all the lines in directions $(\pm 1,\pm 1,\pm 1)$ through integer points.
We now show that these are all the order $3$ rotations in the group,
so that $T=T_1$. For if not, there would be another order $3$ rotation
about a ``new'' axis, and this would differ only by a translation
from a rotation of $T_1$ about a parallel old axis.

However, we show that any translation in the group must be a translation
of the bcc lattice, which preserves the set of old axes.
Our assumptions imply the minimality condition
that if $(t,t,t)$ is a translation in the group,
then $t$ must be a multiple of $\h$.

For since the point group contains the above element $r$, if there is a
translation through $(a,b,c)$, there are others through $(b,c,a)$ and
$(c,a,b)$ and so one through $(a+b+c,a+b+c,a+b+c)$, showing that
$a+b+c$ must be a multiple of $\h$. Similarly, $\pm a\pm  b \pm c$
must be a multiple of $\h$ for all choices of sign, and $2a$, $2b$ and $2c$
must also be multiples of $\h$.

Now transfer the origin to the nearest point of the bcc lattice to $(a,b,c)$
and then change the signs of the coordinates axes, if necessary, to make
$a$, $b$ and $c$ positive.
Up to permutations of the coordinates, $(a,b,c)$ is now
one of
$$
\begin{array}{c}
	(0,0,0),\, (0,0,\h),\, (\q,\q,0),\, (\q,\q,h),\, (\h,\h,0)\,
	\mbox{~or~}\,
	(\h,\h,\h),
\end{array}
$$
the first and last of which are in the bcc lattice and so preserve the set
of old axes.
If any of the other four corresponded to a translation in the group, there
would be a new axis through it in the direction $(1,1,-1)$, but this line
would pass through the point
$$
\begin{array}{c}
	(\q,\q,\q),\, (\e,\e,\e),\, (\E,\E,\E),\, \mbox{~or~}\,  (\q,\q,\q),
\end{array}
$$
respectively, contradicting the minimality condition.

\section*{Appendix III: Details of the Enumeration}

The enumeration proceeds by assigning coupling maps $k\pm$
to the operations in all possible ways that yield homomorphisms
modulo $K$. The following remarks are helpful:

\begin{enumerate}
\item
We often use {\em reduced} names, obtained by omitting
subscripts, when their values are unimportant
or determined.
We have already remarked that the subscript on the digit $A$ between
two punctuation marks is $0$ or $\frac{A}{2}$ according as these are the
same or different.
In the reduced name  $(\s A_a B_b \dots C_c)$
the omitted subscript on $\s$ is necessarily
$\frac{a}{A}+\frac{b}{B}+\dots +\frac{c}{C}$.
\item
In $\up{\lambda}\s\up{P}A\up{Q}B\up{R}\dots C\up{S}$ if $A$ (say)
is odd, then $P$ and $Q$ are conjugate. This entails that the
punctuation marks (if any) on the two sides of $A$ are equal.
So for example, in the case of $\s\up{P}6\up{Q}3\up{R}2$
(where it is $B$ that is odd)
and assuming that all maps  couple to $\p$ elements we get only $4$ cases
$\s\d6\d3\d2$, $\s\d6\dd3\dd2$, $\s\dd6\d3\d2$ and $\s\dd6\dd3\dd2$,
rather than $8$.

\item Often certain parameters can be freely varied without effecting
the truth of the relations.
For example, this happens for $\o\up{X}\up{Y}$: $1=[X,Y]$
when $X\mapsto x\p$ and $Y\mapsto y\p$, since all $k\p$ maps commute.
So this leads to a single case, for which we
choose the couplings $X\mapsto 0\p$ and $Y\mapsto 0\p$ in Table~1.

\item It suffices to work up to symmetry. We note several symmetries:

\begin{enumerate}
\item Changing the sign of the $z$ coordinates: this negates all the constants
in the maps $k\p$ and $k\m$.
This has obvious effects on our names, for example
$(\s 3_1 3_1 3_1)=(\s 3_2 3_2 3_2)$, since
$\t $ and $\T $ are negatives modulo 1.

\item Permuting certain generators: E.g.\ in $ \s\up{P}4\up{Q}4\up{R}2$
we can interchange $P$ and $R$,
in $\s\up{P}3\up{Q}3\up{R}3$ we can apply any permutation, and
in $\s\up{P}2\up{Q}2\up{R}2\up{S}2$ we can cyclically permute $P,Q,R,S$
or reverse their order.
So we have the equalities:
$$
\begin{array}{ccc}
(\s\d4\d4\dd2) & = & (\s\dd4\d4\d2)\\
(\s 3_0 3_1 3_2) & = & (\s 3_1 3_2 3_0)\\
(\s 2_0 2_0 2_1 2_1) & = & (\s 2_0 2_1 2_1 2_0).\\
\end{array}
$$

\item Gyrations can be listed in any order. So in the interval-fiber case
$[2222]$ (where all maps couple to $0\p$ or $\h\p$),
the relation $1=\gamma\delta\epsilon\zeta$ entails that in fact
an even number couple to $\h\p$, leaving just three cases:
$$
\begin{array}{ccl}
0\p 0\p 0\p 0\p & = &
[2_0 2_0 2_0 2_0]\\
0\p 0\p \h\p \h\p &=&
[2_0 2_0 2_1 2_1]\\
\h\p \h\p \h\p \h\p &=&
[2_1 2_1 2_1 2_1].\\
\end{array}
$$

Names that differ only by the obvious symmetries we have described
so far will be regarded as equal.

\item There are more subtle cases: in $\up{\gamma}2\s\up{P}2\up{Q}2\up{R}$,
we eliminated
$R$ as a generator, since $R=\gamma P \gamma ^{-1}$.
But to explore the symmetry it is best to include both $R$ and
$S=\gamma Q\gamma^{-1}$.
We will study the case when all generators couple to $\m$ elements,
say $\gamma\mapsto0\m$, $P\mapsto p\m$, $Q\mapsto q\m$, $R\mapsto r \m$ and
$S\mapsto s\m$.

Then we have $r= -p$, $s= -q$ (modulo 1) and the symmetry permutes
$p,q,r,s$ cyclically, giving the equivalences and reduced names below:

$$
\begin{array}{ccl}
\wrapper{0\m 0\m 0\m 0\m & (2\bar{\s}_0 2_0 2_0)}\\
\wrapper{\h\m \h\m \h\m \h\m &  (2\bar{\s}_1 2_0 2_0)}\\
\multicolumn{2}{l}{\left .\begin{array}{ll}
\h\m 0\m \h\m 0\m &    (2\bar{\s}_0 2_1 2_1)\\
0\m \h\m 0\m \h\m &     (2\bar{\s}_1 2_1 2_1)
\end{array}\right\}}
& (2\bar{\s}2_1 2_1)\\
\multicolumn{2}{l}{\left .\begin{array}{ll}
\q\m \q\m \Q\m \Q\m  &   (2\bar{\s}_0 2_0 2_1)\\
\q\m \Q\m \Q\m \q\m  &   (2\bar{\s}_0 2_1 2_0)\\
\Q\m \Q\m \q\m \q\m  &   (2\bar{\s}_1 2_0 2_1)\\
\Q\m \q\m \q\m \Q\m  &   (2\bar{\s}_1 2_1 2_0)
\end{array}\right\}}
&
(2\bar{\s} 2_0 2_1)
\end{array}
$$

The case $\up{\gamma}4\s\up{P}2\up{Q}$ is similar, with a symmetry
interchanging $P$ and $Q=\gamma P \gamma^{-1}$.
The solution in which the generators $\gamma,P,Q$ all couple to $\m$
elements $0\m$, $p\m$ and $q\m$ are:

$$
\begin{array}{ccc}
\wrapper{0\m 0\m & (4 \bar {\s}_0 2_0)}\\
\wrapper{\h\m \h\m  &  (4 \bar{ \s}_1 2_0)}\\
\multicolumn{2}{l}{\left .
\begin{array}{cc}
\q\m \Q\m  & (4 \bar{ \s}_0 2_1)\\
\Q\m \q\m  & (4 \bar{ \s}_1 2_1)
\end{array}\right \}}
& (4 \bar{ \s} 2_1)\\
\end{array}
$$

\item A still more subtle symmetry arises in the presence of 
a cross cap.

We draw the (reduced) set of generators for
$\up{\gamma}2\up{\delta}2\x \up{Z}$ in Figure~\ref{basechange-fig}(a).
Figure~\ref{basechange-fig}(b) shows that an alternative set of
generators is
$\gamma'=\gamma$, $\delta'=\delta Z\delta^{-1} Z^{-1}\delta^{-1}$
and $Z'=\delta Z$.

\psfragscanon
\begin{figure}[c]
\begin{center}
\begin{tabular}{cc}
\epsfig{file=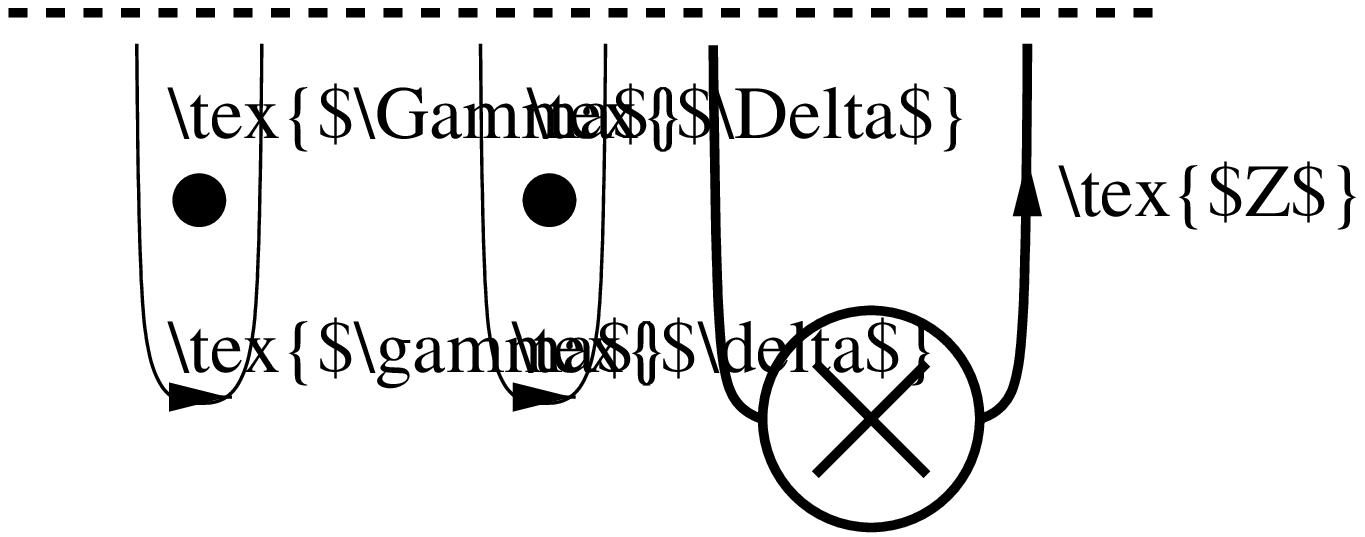,width=7cm}&
\epsfig{file=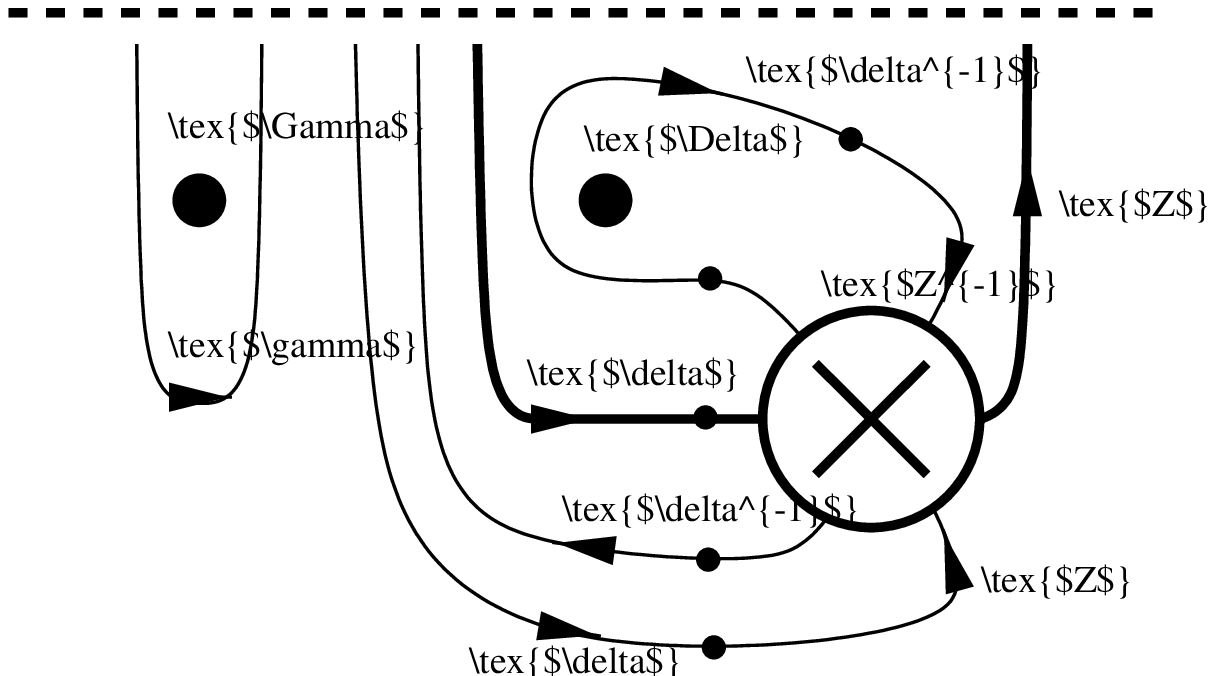,width=7cm}\\
$1=\gamma^2=\delta^2=\gamma\delta Z$&\\
(a)& (b)\\
\end{tabular}
\end{center}
\caption{$\up{\gamma}2\up{\delta}2\x \up{Z}$ possesses a symmetry
that transforms the set of generators indicated in (a) to
the set of generators depicted in (b).}\label{basechange-fig}
\end{figure}
\psfragscanoff

This is because the path $Z$ does not separate the plane,
so the second gyration point $\Delta$ can be pulled through it.
So the paths $\delta Z \delta^{-1}Z^{-1}\delta^{-1}$ and $\delta Z$
of the figure are topologically like $\delta$ and $Z$.

The significant fact here is that $Z$ gets multiplied by $\delta$;
the replacement of $\delta$ by $\delta'$ usually has no effect since
$\delta'$ is conjugate to $\delta^{-1}$.
So in the presence of a $2_1$, $\x_0$ is equivalent to $\x_1$ and in the
presence of a $2$, $\bar{\x}$ is equivalent to $\x$.

This observation reduces the preliminary list of 14 fibrations of $22\x$
to $10$:

\begin{tabular}[t]{ll}
\multicolumn{1}{l}{$\left.\begin{tabular}{l}$[2_0 2_0 \x_0]$\end{tabular}
\right.$} &\\
\multicolumn{1}{l}{$\left.\begin{tabular}{l}$[2_0 2_0 \x_1]$\end{tabular}
\right.$} &\\
\multicolumn{1}{l}{$\left.\begin{tabular}{l}
$[2_1 2_1 \x_0]$\\
$[2_1 2_1 \x_1]$\\ 
\end{tabular}\right\}$
}
&$[2_1 2_1 \x]$ \\
\end{tabular}
\hfill
\begin{tabular}[t]{ll}
\multicolumn{1}{l}{$\left.\begin{tabular}{l}$(2_0 2_0 \x_0)$\end{tabular}
\right.$}&\\
\multicolumn{1}{l}{$\left.\begin{tabular}{l}$(2_0 2_0 \x_1)$\end{tabular}
\right.$}&\\
\multicolumn{1}{l}{$\left.\begin{tabular}{l}
$\left.(2_0 2_1 \x_0)\right.$\\
$\left.(2_0 2_1 \x_1)\right.$\\
\end{tabular}\right\}$ }
& $(2_0 2_1 \x)$\\
\multicolumn{1}{l}{$\left.\begin{tabular}{l}
$(2_1 2_1 \x_0)$\\
$(2_1 2_1 \x_1)$\\
\end{tabular}\right\}$ }
& $(2_1 2_1 \x)$\\
\end{tabular}\hfill
\begin{tabular}[t]{ll}
\multicolumn{1}{l}{$\left.\begin{tabular}{l}$(2_0 2_0 \bar{\x})$\end{tabular}
\right.$}&\\
\multicolumn{1}{l}{$\left.\begin{tabular}{l}$(2_1 2_1 \bar{\x})$\end{tabular}
\right.$}&\\
&\\
\multicolumn{1}{l}{$\left.\begin{tabular}{l}
$(2 2 \x)$\\
$(2 2 \bar{\x})$\\
\end{tabular}\right\}$ }
& $(2 2 \x)$\\
\end{tabular}\\
It turns out that these are all distinct.
\end{enumerate}

\end{enumerate}

\newpage
\bibliographystyle{alpha}
\bibliography{maths}

\newpage

\section*{Captions for main tables}

\paragraph{Table 1}
This table consists of 17 blocks.
Each is headed by a plane group $H$ and a set of relations for $H$.
This if followed by a number of lines corresponding to different fibrations
over $H$.
In each such line, we list the fibrifold name for the resulting space group $G$
(first column), the appropriate couplings for the generators of
$H$ (second column), the point group of $G$ (third column) and
finally, the IT number of $G$ (fourth column).
In all three tables, ($\ddagger$) indicates that a group has two
enantiomorphous forms and we report only the lower IT number.

\paragraph{Table 2a}
Here we list the 35 irreducible groups, organized by their point groups.
Each line describes one such group $G$,
listing its primary name (column one), its international number and  name
(column two) and its secondary name(s) (column three).

\paragraph{Table 2b}
Here we list the 184 reducible groups, organized by their point groups.
Each line describes one such group $G$,
listing its primary name (column one),
its international number and  name (column two)
and, its secondary name(s), if it has any (column three).
\textwidth18cm
\hoffset=-3cm
\textheight25cm
\voffset=-3cm

%
%
%
%
%
%
%
%
%
\twocolumn
\begin{center}

\tablefirsthead{\hline \multicolumn{4}{|c|}{\bf Table 1 (\today)}\\\hline}

\tablehead{\hline \multicolumn{4}{|l|}{\em Table 1 continued}\\\hline}

\tabletail{\hline}

\begin{supertabular}{|c|l|c|c|}

\subheader
{$\s 632$}
{$\s \up{P}6\up{Q}3\up{R}2$}
{$1=P^2=(PQ)^6=Q^2=(QR)^3=R^2=(RP)^2$}
{$P$ {\hspace{3pt}} $Q$ {\hspace{3pt}} $R$ {\hspace{3pt}} \id}

$[\s \d6\d3\d2]$ & $\xo\p \xo\p \xo\p \xR $ & $\s 226$ & $191$\\
$[\s \dd6\d3\d2]$ & $\xh\p \xo\p \xo\p \xR $ & $\s 226$ & $194$\\
$[\s \d6\dd3\dd2]$ & $\xo\p \xh\p \xh\p \xR $ & $\s 226$ & $193$\\
$[\s \dd6\dd3\dd2]$ & $\xh\p \xh\p \xh\p \xR $ & $\s 226$ & $192$\\
&&&\\

$(\s \d6\d3\d2)$ & $\xo\p \xo\p \xo\p $ & $\s 66$ & $183$\\
$(\s \dd6\d3\d2)$ & $\xh\p \xo\p \xo\p $ & $\s 66$ & $186$\\
$(\s \d6\dd3\dd2)$ & $\xo\p \xh\p \xh\p $ & $\s 66$ & $185$\\
$(\s \dd6\dd3\dd2)$ & $\xh\p \xh\p \xh\p $ & $\s 66$ & $184$\\
&&&\\

$(\s  6\d3\d2)$ & $\xo\m \xo\p \xo\p $ & $2\s 3$ & $164$\\
$(\s  6\dd3\dd2)$ & $\xo\m \xh\p \xh\p $ & $2\s 3$ & $165$\\
&&&\\

$(\s \d6\n 3_0 2)$ & $\xo\p \xo\m \xo\m $ & $2\s 3$ & $162$\\
$(\s \d6\n 3_1 2)$ & $\xo\p \xt\m \xo\m $ & $2\s 3$ & $166$\\
$(\s \dd6\n 3_0 2)$ & $\xh\p \xo\m \xo\m $ & $2\s 3$ & $163$\\
$(\s \dd6\n 3_1 2)$ & $\xh\p \xt\m \xo\m $ & $2\s 3$ & $167$\\
&&&\\

$(\s  6_0 3_0 2_0)$ & $\xo\m \xo\m \xo\m $ & $226$ & $177$\\
$(\s  6_1 3_1 2_1)$ & $\xh\m \xt\m \xo\m $ & $226$ & $178$\\
$(\s  6_2 3_2 2_0)$ & $\xo\m \xT\m \xo\m $ & $226$ & $180$\\
$(\s  6_3 3_0 2_1)$ & $\xh\m \xo\m \xo\m $ & $226$ & $182$\\
&&&\\

\subheader
{$632$}
{$\up{\gamma}6\up{\delta}3\up{\epsilon}2$}
{$1=\gamma^6=\delta^3=\epsilon^2=\gamma\delta\epsilon$}
{$\gamma$ \, $\delta$ \,\, $\epsilon$ \, \id}

$[6_0 3_0 2_0]$ & $\xo\p \xo\p \xo\p \xR $ & $6\s $ & $175$\\
$[6_3 3_0 2_1]$ & $\xh\p \xo\p \xh\p \xR $ & $6\s $ & $176$\\
&&&\\

$(6_0 3_0 2_0)$ & $\xo\p \xo\p \xo\p $ & $66$ & $168$\\
$(6_1 3_1 2_1)$ & $\xs\p \xt\p \xh\p $ & $66$ & $169$\\
$(6_2 3_2 2_0)$ & $\xt\p \xT\p \xo\p $ & $66$ & $171$\\
$(6_3 3_0 2_1)$ & $\xh\p \xo\p \xh\p $ & $66$ & $173$\\
&&&\\

$(6\n 3_0 2)$ & $\xo\m \xo\p \xo\m $ & $3\x$ & $147$\\
$(6\n 3_1 2)$ & $\xt\m \xt\p \xo\m $ & $3\x$ & $148$\\
&&&\\
\end{supertabular}

\tablefirsthead{\hline \multicolumn{4}{|l|}{\em Table 1 continued}\\\hline}
\begin{supertabular}{|c|l|c|c|}

\subheader
{$\s 442$}
{$\s \up{P}4\up{Q}4\up{R}2$}
{$1=P^2=(PQ)^4=Q^2=(QR)^4=R^2=(RP)^2$}
{$P$ {\hspace{3pt}} $Q$ {\hspace{3pt}} $R$ {\hspace{3pt}} \id}

$[\s \d4\d4\d2]$ & $\xo\p \xo\p \xo\p \xR $ & $\s 224$ & $123$\\
$[\s \d4\d4\dd2]$ & $\xo\p \xo\p \xh\p \xR $ & $\s 224$ & $139$\\
$[\s \d4\dd4\d2]$ & $\xo\p \xh\p \xo\p \xR $ & $\s 224$ & $131$\\
$[\s \d4\dd4\dd2]$ & $\xo\p \xh\p \xh\p \xR $ & $\s 224$ & $140$\\
$[\s \dd4\d4\dd2]$ & $\xh\p \xo\p \xh\p \xR $ & $\s 224$ & $132$\\
$[\s \dd4\dd4\dd2]$ & $\xh\p \xh\p \xh\p \xR $ & $\s 224$ & $124$\\
&&&\\

$(\s \d4\d4\d2)$ & $\xo\p \xo\p \xo\p $ & $\s 44$ & $99$\\
$(\s \d4\d4\dd2)$ & $\xo\p \xo\p \xh\p $ & $\s 44$ & $107$\\
$(\s \d4\dd4\d2)$ & $\xo\p \xh\p \xo\p $ & $\s 44$ & $105$\\
$(\s \d4\dd4\dd2)$ & $\xo\p \xh\p \xh\p $ & $\s 44$ & $108$\\
$(\s \dd4\d4\dd2)$ & $\xh\p \xo\p \xh\p $ & $\s 44$ & $101$\\
$(\s \dd4\dd4\dd2)$ & $\xh\p \xh\p \xh\p $ & $\s 44$ & $103$\\
&&&\\

$(\s  4\d4\d2)$ & $\xo\m \xo\p \xo\p     $ & $\s 224$ & $129$\\
$(\s  4\d4\dd2)$ & $\xo\m \xo\p \xh\p $ & $\s 224$ & $137$\\
$(\s  4\dd4\d2)$ & $\xo\m \xh\p \xo\p $ & $\s 224$ & $138$\\
$(\s  4\dd4\dd2)$ & $\xo\m \xh\p \xh\p $ & $\s 224$ & $130$\\
&&&\\

$(\s \d4\n 4\d2)$ & $\xo\p \xo\m \xo\p $ & $2\s 2$ & $115$\\
$(\s \d4\n 4\dd2)$ & $\xo\p \xo\m \xh\p $ & $2\s 2$ & $121$\\
$(\s \dd4\n 4\dd2)$ & $\xh\p \xo\m \xh\p $ & $2\s 2$ & $116$\\
&&&\\

$(\s  4_0 4\d2)$ & $\xo\m \xo\m \xo\p $ & $\s 224$ & $125$\\
$(\s  4_1 4\d2)$ & $\xq\m \xo\m \xo\p $ & $\s 224$ & $141$\\
$(\s  4_2 4\d2)$ & $\xh\m \xo\m \xo\p $ & $\s 224$ & $134$\\
$(\s  4_0 4\dd2)$ & $\xo\m \xo\m \xh\p $ & $\s 224$ & $126$\\
$(\s  4_1 4\dd2)$ & $\xq\m \xo\m \xh\p $ & $\s 224$ & $142$\\
$(\s  4_2 4\dd2)$ & $\xh\m \xo\m \xh\p $ & $\s 224$ & $133$\\
&&&\\

$(\s  4\d4\n 2_0)$ & $\xo\m \xo\p \xo\m $ & $2\s 2$ & $111$\\
$(\s  4\d4\n 2_1)$ & $\xh\m \xo\p \xo\m $ & $2\s 2$ & $119$\\
$(\s  4\dd4\n 2_0)$ & $\xo\m \xh\p \xo\m $ & $2\s 2$ & $112$\\
$(\s  4\dd4\n 2_1)$ & $\xh\m \xh\p \xo\m $ & $2\s 2$ & $120$\\
&&&\\

$(\s  4_0 4_0 2_0)$ & $\xo\m \xo\m \xo\m $ & $224$ & $89$\\
$(\s  4_1 4_1 2_1)$ & $\xh\m \xq\m \xo\m $ & $224$ & $91$\\
$(\s  4_2 4_2 2_0)$ & $\xo\m \xh\m \xo\m $ & $224$ & $93$\\
$(\s  4_2 4_0 2_1)$ & $\xh\m \xo\m \xo\m $ & $224$ & $97$\\
$(\s  4_3 4_1 2_0)$ & $\xo\m \xq\m \xo\m $ & $224$ & $98$\\
&&&\\
\end{supertabular}

\begin{supertabular}{|c|l|c|c|}

\subheader
{$4\s2$}
{$ \up{\gamma}4\s\up{P}2$}
{$1=\gamma^4=P^2=[P,\gamma]^2$}
{$\gamma$ \, $P$  \id}

$[4_0 \s \d2]$ & $\xo\p \xo\p \xR $ & $\s 224$ & $127$\\
$[4_2 \s \d2]$ & $\xh\p \xo\p \xR $ & $\s 224$ & $136$\\
$[4_0 \s \dd2]$ & $\xo\p \xh\p \xR $ & $\s 224$ & $128$\\
$[4_2 \s \dd2]$ & $\xh\p \xh\p \xR $ & $\s 224$ & $135$\\
&&&\\

$(4_0 \s \d2)$ & $\xo\p \xo\p $ & $\s 44$ & $100$\\
$(4_1 \s \d2)$ & $\xq\p \xo\p $ & $\s 44$ & $109$\\
$(4_2 \s \d2)$ & $\xh\p \xo\p $ & $\s 44$ & $102$\\
$(4_0 \s \dd2)$ & $\xo\p \xh\p $ & $\s 44$ & $104$\\
$(4_1 \s \dd2)$ & $\xq\p \xh\p $ & $\s 44$ & $110$\\
$(4_2 \s \dd2)$ & $\xh\p \xh\p $ & $\s 44$ & $106$\\
&&&\\

$(4 \bar{\s} \d2)$ & $\xo\m \xo\p $ & $2\s 2$ & $113$\\
$(4 \bar{\s} \dd2)$ & $\xo\m \xh\p $ & $2\s 2$ & $114$\\
&&&\\

$(4_0 \s  2_0)$ & $\xo\p \xo\m $ & $224$ & $90$\\
$(4_1 \s  2_1)$ & $\xq\p \xo\m $ & $224$ & $92$\\
$(4_2 \s  2_0)$ & $\xh\p \xo\m $ & $224$ & $94$\\
&&&\\

$(4 \bar{\s}_0 2_0)$ & $\xo\m \xo\m $ & $2\s 2$ & $117$\\
$(4 \bar{\s}_1 2_0)$ & $\xo\m \xh\m $ & $2\s 2$ & $118$\\
$(4 \bar{\s}  2_1)$ & $\xo\m \xq\m $ & $2\s 2$ & $122$\\
&&&\\

\subheader
{$442$}
{$ \up{\gamma}4\up{\delta}4\up{\epsilon}2$}
{$1=\gamma^4=\delta^4=\epsilon^2=\gamma\delta\epsilon$}
{$\gamma$ \,\, $\delta$ \, \,$\epsilon$ \,\id}

$[4_0 4_0 2_0]$ & $\xo\p \xo\p \xo\p \xR $ & $4\s $ & $83$\\
$[4_2 4_2 2_0]$ & $\xh\p \xh\p \xo\p \xR $ & $4\s $ & $84$\\
$[4_2 4_0 2_1]$ & $\xh\p \xo\p \xh\p \xR $ & $4\s $ & $87$\\
&&&\\

$(4_0 4_0 2_0)$ & $\xo\p \xo\p \xo\p $ & $44$ & $75$\\
$(4_1 4_1 2_1)$ & $\xq\p \xq\p \xh\p $ & $44$ & $76$\\
$(4_2 4_2 2_0)$ & $\xh\p \xh\p \xo\p $ & $44$ & $77$\\
$(4_2 4_0 2_1)$ & $\xh\p \xo\p \xh\p $ & $44$ & $79$\\
$(4_3 4_1 2_0)$ & $\xQ\p \xq\p \xo\p $ & $44$ & $80$\\
&&&\\

$(4\n 4_0 2)$ & $\xo\p \xo\m \xo\m $ & $4\s $ & $85$\\
$(4\n 4_1 2)$ & $\xq\p \xo\m \xq\m $ & $4\s $ & $88$\\
$(4\n 4_2 2)$ & $\xh\p \xo\m \xh\m $ & $4\s $ & $86$\\
&&&\\

$(4\n 4\n 2_0)$ & $\xo\m \xo\m \xo\p $ & $2\x$ & $81$\\
$(4\n 4\n 2_1)$ & $\xh\m \xo\m \xh\p $ & $2\x$ & $82$\\
&&&\\
\end{supertabular}

\begin{supertabular}{|c|l|c|c|}

\subheader
{$\s333$}
{$ \s\up{P}3\up{Q}3\up{R}3$}
{$1=P^2=(PQ)^3=Q^2=(QR)^3=R^2=(RP)^3$}
{$P$ {\hspace{3pt}} $Q$ {\hspace{3pt}} $R$ {\hspace{3pt}} \id}

$[\s \d3\d3\d3]$ & $\xo\p \xo\p \xo\p \xR $ & $\s 223$ & $187$\\
$[\s \dd3\dd3\dd3]$ & $\xh\p \xh\p \xh\p \xR $ & $\s 223$ & $188$\\
&&&\\

$(\s \d3\d3\d3)$ & $\xo\p \xo\p \xo\p $ & $\s 33$ & $156$\\
$(\s \dd3\dd3\dd3)$ & $\xh\p \xh\p \xh\p $ & $\s 33$ & $158$\\
&&&\\

$(\s  3_0 3_0 3_0)$ & $\xo\m \xo\m \xo\m $ & $223$ & $149$\\
$(\s  3_1 3_1 3_1)$ & $\xT\m \xt\m \xo\m $ & $223$ & $151$\\
$(\s  3_0 3_1 3_2)$ & $\xt\m \xt\m \xo\m $ & $223$ & $155$\\
&&&\\

\subheader
{$3\s3$}
{$ \up{\gamma}3\s\up{P}3$}
{$1=\gamma^3=P^2=[P,\gamma]^3$}
{$\gamma$ \, $P$  \id}

$[3_0 \s \d3]$ & $\xo\p \xo\p \xR $ & $\s 223$ & $189$\\
$[3_0 \s \dd3]$ & $\xo\p \xh\p \xR $ & $\s 223$ & $190$\\
&&&\\

$(3_0 \s \d3)$ & $\xo\p \xo\p $ & $\s 33$ & $157$\\
$(3_1 \s \d3)$ & $\xt\p \xo\p $ & $\s 33$ & $160$\\
$(3_0 \s \dd3)$ & $\xo\p \xh\p $ & $\s 33$ & $159$\\
$(3_1 \s \dd3)$ & $\xt\p \xh\p $ & $\s 33$ & $161$\\
&&&\\

$(3_0 \s  3_0)$ & $\xo\p \xo\m $ & $223$ & $150$\\
$(3_1 \s  3_1)$ & $\xt\p \xo\m $ & $223$ & $152$\\
&&&\\

\subheader
{$333$}
{$\up{\gamma}3\up{\delta}3\up{\epsilon}3$}
{$1=\gamma^3=\delta^3=\epsilon^3=\gamma\delta\epsilon$}
{$\gamma$ \, $\delta$ \,\, $\epsilon$ \, \id}

$[3_0 3_0 3_0]$ & $\xo\p \xo\p \xo\p \xR $ & $3\s $ & $174$\\
&&&\\

$(3_0 3_0 3_0)$ & $\xo\p \xo\p \xo\p $ & $33$ & $143$\\
$(3_1 3_1 3_1)$ & $\xt\p \xt\p \xt\p $ & $33$ & $144$\\
$(3_0 3_1 3_2)$ & $\xo\p \xt\p \xT\p $ & $33$ & $146$\\
&&&\\
\end{supertabular}

\begin{supertabular}{|c|l|c|c|}

\subheader
{$\s2222$}
{$ \s\up{P}2\up{Q}2\up{R}2\up{S}2$}
{$1=P^2=(PQ)^2=Q^2=(QR)^2=R^2=(RS)^2=S^2=(SP)^2$}
{$P$ {\hspace{3pt}} $Q$ {\hspace{3pt}} $R$ {\hspace{3pt}} $S$ {\hspace{3pt}} \id}

$[\s \d2\d2\d2\d2]$ & $\xo\p \xo\p \xo\p \xo\p \xR $ & $\s 222$ & $47$\\
$[\s \d2\d2\d2\dd2]$ & $\xo\p \xo\p \xo\p \xh\p \xR $ & $\s 222$ & $65$\\
$[\s \d2\d2\dd2\dd2]$ & $\xo\p \xo\p \xh\p \xh\p \xR $ & $\s 222$ & $69$\\
$[\s \d2\dd2\d2\dd2]$ & $\xo\p \xh\p \xo\p \xh\p \xR $ & $\s 222$ & $51$\\
$[\s \d2\dd2\dd2\dd2]$ & $\xo\p \xh\p \xh\p \xh\p \xR $ & $\s 222$ & $67$\\
$[\s \dd2\dd2\dd2\dd2]$ & $\xh\p \xh\p \xh\p \xh\p \xR $ & $\s 222$ & $49$\\
&&&\\

$(\s \d2\d2\d2\d2)$ & $\xo\p \xo\p \xo\p \xo\p $ & $\s 22$ & $25$\\
$(\s \d2\d2\d2\dd2)$ & $\xo\p \xo\p \xo\p \xh\p $ & $\s 22$ & $38$\\
$(\s \d2\d2\dd2\dd2)$ & $\xo\p \xo\p \xh\p \xh\p $ & $\s 22$ & $42$\\
$(\s \d2\dd2\d2\dd2)$ & $\xo\p \xh\p \xo\p \xh\p $ & $\s 22$ & $26$\\
$(\s \d2\dd2\dd2\dd2)$ & $\xo\p \xh\p \xh\p \xh\p $ & $\s 22$ & $39$\\
$(\s \dd2\dd2\dd2\dd2)$ & $\xh\p \xh\p \xh\p \xh\p $ & $\s 22$ & $27$\\
&&&\\

$(\s  2\d2\d2\d2)$ & $\xo\p \xo\p \xo\p \xo\m $ & $\s 222$ & $51$\\
$(\s  2\d2\d2\dd2)$ & $\xo\p \xo\p \xh\p \xo\m $ & $\s 222$ & $63$\\
$(\s  2\d2\dd2\d2)$ & $\xo\p \xh\p \xo\p \xo\m $ & $\s 222$ & $55$\\
$(\s  2\d2\dd2\dd2)$ & $\xo\p \xh\p \xh\p \xo\m $ & $\s 222$ & $64$\\
$(\s  2\dd2\d2\dd2)$ & $\xh\p \xo\p \xh\p \xo\m $ & $\s 222$ & $57$\\
$(\s  2\dd2\dd2\dd2)$ & $\xh\p \xh\p \xh\p \xo\m $ & $\s 222$ & $54$\\
&&&\\

$(\s  2_0 2\d2\d2)$ & $\xo\m \xo\m \xo\p \xo\p $ & $\s 222$ & $67$\\
$(\s  2_1 2\d2\d2)$ & $\xh\m \xo\m \xo\p \xo\p $ & $\s 222$ & $74$\\
$(\s  2_0 2\d2\dd2)$ & $\xo\m \xo\m \xo\p \xh\p $ & $\s 222$ & $72$\\
$(\s  2_1 2\d2\dd2)$ & $\xh\m \xo\m \xo\p \xh\p $ & $\s 222$ & $64$\\
$(\s  2_0 2\dd2\dd2)$ & $\xo\m \xo\m \xh\p \xh\p $ & $\s 222$ & $68$\\
$(\s  2_1 2\dd2\dd2)$ & $\xh\m \xo\m \xh\p \xh\p $ & $\s 222$ & $73$\\
&&&\\

$(\s  2\d2\n 2\d2)$ & $\xo\m \xo\p \xo\m \xo\p $ & $2\s $ & $10$\\
$(\s  2\d2\n 2\dd2)$ & $\xo\m \xo\p \xo\m \xh\p $ & $2\s $ & $12$\\
$(\s  2\dd2\n 2\dd2)$ & $\xo\m \xh\p \xo\m \xh\p $ & $2\s $ & $13$\\
&&&\\

$(\s  2_0 2_0 2\d2)$ & $\xo\m \xo\m \xo\m \xo\p $ & $\s 222$ & $49$\\
$(\s  2_0 2_1 2\d2)$ & $\xo\m \xo\m \xh\m \xo\p $ & $\s 222$ & $66$\\
$(\s  2_1 2_1 2\d2)$ & $\xh\m \xo\m \xh\m \xo\p $ & $\s 222$ & $53$\\
$(\s  2_0 2_0 2\dd2)$ & $\xo\m \xo\m \xo\m \xh\p $ & $\s 222$ & $50$\\
$(\s  2_0 2_1 2\dd2)$ & $\xo\m \xo\m \xh\m \xh\p $ & $\s 222$ & $68$\\
$(\s  2_1 2_1 2\dd2)$ & $\xh\m \xo\m \xh\m \xh\p $ & $\s 222$ & $54$\\
&&&\\

$(\s  2_0 2_0 2_0 2_0)$ & $\xo\m \xo\m \xo\m \xo\m $ & $222$ & $16$\\
$(\s  2_0 2_0 2_1 2_1)$ & $\xo\m \xo\m \xo\m \xh\m $ & $222$ & $21$\\
$(\s  2_0 2_1 2_0 2_1)$ & $\xo\m \xo\m \xh\m \xh\m $ & $222$ & $22$\\
$(\s  2_1 2_1 2_1 2_1)$ & $\xo\m \xh\m \xo\m \xh\m $ & $222$ & $17$\\
&&&\\
\end{supertabular}

\begin{supertabular}{|c|l|c|c|}

\subheader
{$2\s22$}
{$ \up{\gamma}2\s\up{P}2\up{Q}2$}
{$1=\gamma^2=P^2=(PQ)^2=Q^2=(Q\gamma P\gamma^{-1})^2$}
{$\gamma$ \, $P$ {\hspace{3pt}} $Q$ {\hspace{3pt}} \id}

$[2_0 \s \d2\d2]$ & $\xo\p \xo\p \xo\p \xR $ & $\s 222$ & $65$\\
$[2_1 \s \d2\d2]$ & $\xh\p \xo\p \xo\p \xR $ & $\s 222$ & $71$\\
$[2_0 \s \d2\dd2]$ & $\xo\p \xo\p \xh\p \xR $ & $\s 222$ & $74$\\
$[2_1 \s \d2\dd2]$ & $\xh\p \xo\p \xh\p \xR $ & $\s 222$ & $63$\\
$[2_0 \s \dd2\dd2]$ & $\xo\p \xh\p \xh\p \xR $ & $\s 222$ & $66$\\
$[2_1 \s \dd2\dd2]$ & $\xh\p \xh\p \xh\p \xR $ & $\s 222$ & $72$\\
&&&\\

$(2_0 \s \d2\d2)$ & $\xo\p \xo\p \xo\p $ & $\s 22$ & $35$\\
$(2_1 \s \d2\d2)$ & $\xh\p \xo\p \xo\p $ & $\s 22$ & $44$\\
$(2_0 \s \d2\dd2)$ & $\xo\p \xo\p \xh\p $ & $\s 22$ & $46$\\
$(2_1 \s \d2\dd2)$ & $\xh\p \xo\p \xh\p $ & $\s 22$ & $36$\\
$(2_0 \s \dd2\dd2)$ & $\xo\p \xh\p \xh\p $ & $\s 22$ & $37$\\
$(2_1 \s \dd2\dd2)$ & $\xh\p \xh\p \xh\p $ & $\s 22$ & $45$\\
&&&\\

$(2_0 \s  2\d2)$ & $\xo\p \xo\p \xo\m $ & $\s 222$ & $53$\\
$(2_1 \s  2\d2)$ & $\xh\p \xo\p \xo\m $ & $\s 222$ & $58$\\
$(2_0 \s  2\dd2)$ & $\xo\p \xh\p \xo\m $ & $\s 222$ & $52$\\
$(2_1 \s  2\dd2)$ & $\xh\p \xh\p \xo\m $ & $\s 222$ & $60$\\
&&&\\

$(2 \bar{\s} \d2\d2)$ & $\xo\m \xo\p \xo\p $ & $\s 222$ & $59$\\
$(2 \bar{\s} \d2\dd2)$ & $\xo\m \xo\p \xh\p $ & $\s 222$ & $62$\\
$(2 \bar{\s} \dd2\dd2)$ & $\xo\m \xh\p \xh\p $ & $\s 222$ & $56$\\
&&&\\

$(2 \bar{\s}  2\d2)$ & $\xo\m \xo\m \xo\p $ & $2\s $ & $12$\\
$(2 \bar{\s}  2\dd2)$ & $\xo\m \xo\m \xh\p $ & $2\s $ & $15$\\
&&&\\

$(2_0 \s  2_0 2_0)$ & $\xo\p \xo\m \xo\m $ & $222$ & $21$\\
$(2_1 \s  2_0 2_0)$ & $\xh\p \xo\m \xo\m $ & $222$ & $23$\\
$(2_0 \s  2_1 2_1)$ & $\xo\p \xh\m \xo\m $ & $222$ & $24$\\
$(2_1 \s  2_1 2_1)$ & $\xh\p \xh\m \xo\m $ & $222$ & $20$\\
&&&\\

$(2 \bar{\s}_0 2_0 2_0)$ & $\xo\m \xo\m \xo\m $ & $\s 222$ & $50$\\
$(2 \bar{\s}_1 2_0 2_0)$ & $\xo\m \xh\m \xh\m $ & $\s 222$ & $48$\\
$(2 \bar{\s}  2_0 2_1)$ & $\xo\m \xq\m \xq\m $ & $\s 222$ & $70$\\
$(2 \bar{\s}  2_1 2_1)$ & $\xo\m \xo\m \xh\m $ & $\s 222$ & $52$\\
&&&\\
\end{supertabular}

\begin{supertabular}{|c|l|c|c|}

\subheader
{$22\s$}
{$ \up{\gamma}2\up{\delta}2\s\up{P}$}
{$1=\gamma^2=\delta^2=P^2=[P,\gamma\delta]$}
{$\gamma$ \, $\delta$ \, $P$ {\hspace{3pt}} \id}

$[2_0 2_0 \s \d]$ & $\xo\p \xo\p \xo\p \xR $ & $\s 222$ & $51$\\
$[2_0 2_1 \s \d]$ & $\xo\p \xh\p \xo\p \xR $ & $\s 222$ & $63$\\
$[2_1 2_1 \s \d]$ & $\xh\p \xh\p \xo\p \xR $ & $\s 222$ & $59$\\
$[2_0 2_0 \s \dd]$ & $\xo\p \xo\p \xh\p \xR $ & $\s 222$ & $53$\\
$[2_0 2_1 \s \dd]$ & $\xo\p \xh\p \xh\p \xR $ & $\s 222$ & $64$\\
$[2_1 2_1 \s \dd]$ & $\xh\p \xh\p \xh\p \xR $ & $\s 222$ & $57$\\
&&&\\

$(2_0 2_0 \s \d)$ & $\xo\p \xo\p \xo\p $ & $\s 22$ & $28$\\
$(2_0 2_1 \s \d)$ & $\xo\p \xh\p \xo\p $ & $\s 22$ & $40$\\
$(2_1 2_1 \s \d)$ & $\xh\p \xh\p \xo\p $ & $\s 22$ & $31$\\
$(2_0 2_0 \s \dd)$ & $\xo\p \xo\p \xh\p $ & $\s 22$ & $30$\\
$(2_0 2_1 \s \dd)$ & $\xo\p \xh\p \xh\p $ & $\s 22$ & $41$\\
$(2_1 2_1 \s \dd)$ & $\xh\p \xh\p \xh\p $ & $\s 22$ & $29$\\
&&&\\

$(2_0 2_0 \s )$ & $\xo\p \xo\p \xo\m $ & $222$ & $17$\\
$(2_0 2_1 \s )$ & $\xo\p \xh\p \xo\m $ & $222$ & $20$\\
$(2_1 2_1 \s )$ & $\xh\p \xh\p \xo\m $ & $222$ & $18$\\
&&&\\

$(2_0 2 \bar{\s} \d)$ & $\xo\p \xo\m \xo\p $ & $\s 222$ & $57$\\
$(2_1 2 \bar{\s} \d)$ & $\xh\p \xo\m \xo\p $ & $\s 222$ & $62$\\
$(2_0 2 \bar{\s} \dd)$ & $\xo\p \xo\m \xh\p $ & $\s 222$ & $60$\\
$(2_1 2 \bar{\s} \dd)$ & $\xh\p \xo\m \xh\p $ & $\s 222$ & $61$\\
&&&\\

$(2_0 2 \bar{\s}_0)$ & $\xo\p \xo\m \xo\m $ & $\s 222$ & $54$\\
$(2_0 2 \bar{\s}_1)$ & $\xo\p \xh\m \xo\m $ & $\s 222$ & $52$\\
$(2_1 2 \bar{\s}_0)$ & $\xh\p \xh\m \xo\m $ & $\s 222$ & $56$\\
$(2_1 2 \bar{\s}_1)$ & $\xh\p \xo\m \xo\m $ & $\s 222$ & $60$\\
&&&\\

$(2\n 2 \s \d)$ & $\xo\m \xo\m \xo\p $ & $2\s $ & $11$\\
$(2\n 2 \s \dd)$ & $\xo\m \xo\m \xh\p $ & $2\s $ & $14$\\
&&&\\

$(2\n 2 \s_0)$ & $\xo\m \xo\m \xo\m $ & $2\s $ & $13$\\
$(2\n 2 \s_1)$ & $\xo\m \xh\m \xo\m $ & $2\s $ & $15$\\
&&&\\
\end{supertabular}

\begin{supertabular}{|c|l|c|c|}

\subheader
{$22\x$}
{$ \up{\gamma}2\up{\delta}2\x\up{Z}$}
{$1=\gamma^2=\delta^2=\gamma\delta Z^2$}
{$\gamma$ \, $\delta$ \, $Z$ {\hspace{3pt}} \id}

$[2_0 2_0 \x_0]$ & $\xo\p \xo\p \xo\p \xR $ & $\s 222$ & $55$\\
$[2_0 2_0 \x_1]$ & $\xo\p \xo\p \xh\p \xR $ & $\s 222$ & $58$\\
$[2_1 2_1 \x]$ & $\xh\p \xh\p \xo\p \xR $ & $\s 222$ & $62$\\
&&&\\

$(2_0 2_0 \x_0)$ & $\xo\p \xo\p \xo\p $ & $\s 22$ & $32$\\
$(2_0 2_0 \x_1)$ & $\xo\p \xo\p \xh\p $ & $\s 22$ & $34$\\
$(2_0 2_1 \x)$ & $\xo\p \xh\p \xq\p $ & $\s 22$ & $43$\\
$(2_1 2_1 \x)$ & $\xh\p \xh\p \xo\p $ & $\s 22$ & $33$\\
&&&\\

$(2_0 2_0 \bar{\x})$ & $\xo\p \xo\p \xo\m $ & $222$ & $18$\\
$(2_1 2_1 \bar{\x})$ & $\xh\p \xh\p \xo\m $ & $222$ & $19$\\
&&&\\

$(2\n 2 \x)$ & $\xo\m \xo\m \xo\p $ & $2\s $ & $14$\\
&&&\\

\subheader
{$2222$}
{$\up{\gamma}2\up{\delta}2\up{\epsilon}2\up{\zeta}2$}
{$1=\gamma^2=\delta^2=\epsilon^2=\zeta^2=\gamma\delta\epsilon\zeta$}
{$\gamma$ \, $\delta$ \,\, $\epsilon$ \,\, $\zeta$ {\hspace{3pt}} \id}

$[2_0 2_0 2_0 2_0]$ & $\xo\p \xo\p \xo\p \xo\p \xR $ & $2\s $ & $10$\\
$[2_0 2_0 2_1 2_1]$ & $\xo\p \xo\p \xh\p \xh\p \xR $ & $2\s $ & $12$\\
$[2_1 2_1 2_1 2_1]$ & $\xh\p \xh\p \xh\p \xh\p \xR $ & $2\s $ & $11$\\
&&&\\

$(2_0 2_0 2_0 2_0)$ & $\xo\p \xo\p \xo\p \xo\p $ & $22$ & $3$\\
$(2_0 2_0 2_1 2_1)$ & $\xo\p \xo\p \xh\p \xh\p $ & $22$ & $5$\\
$(2_1 2_1 2_1 2_1)$ & $\xh\p \xh\p \xh\p \xh\p $ & $22$ & $4$\\
&&&\\

$(2_0 2_0 2\n 2)$ & $\xo\p \xo\p \xo\m \xo\m $ & $2\s $ & $13$\\
$(2_0 2_1 2\n 2)$ & $\xo\p \xh\p \xo\m \xh\m $ & $2\s $ & $15$\\
$(2_1 2_1 2\n 2)$ & $\xh\p \xh\p \xo\m \xo\m $ & $2\s $ & $14$\\
&&&\\

$(2\n 2\n 2\n 2)$ & $\xo\m \xo\m \xo\m \xo\m $ & $\x$ & $2$\\
&&&\\
\end{supertabular}

\begin{supertabular}{|c|l|c|c|}

\subheader
{$\s\s$}
{$\up{\lambda}\s\up{P}\s\up{Q}$}
{$1=P^2=Q^2=[\lambda,P]=[\lambda,Q]$}
{$\lambda$ \, $P$ {\hspace{3pt}} $Q$ {\hspace{3pt}} \id}

$[\s_0 \d\s_0 \d]$ & $\xo\p \xo\p \xo\p \xR $ & $\s 22$ & $25$\\
$[\s_1 \d\s_1 \d]$ & $\xh\p \xo\p \xo\p \xR $ & $\s 22$ & $38$\\
$[\s_0 \d\s_0 \dd]$ & $\xo\p \xo\p \xh\p \xR $ & $\s 22$ & $35$\\
$[\s_1 \d\s_1 \dd]$ & $\xh\p \xo\p \xh\p \xR $ & $\s 22$ & $42$\\
$[\s_0 \dd\s_0 \dd]$ & $\xo\p \xh\p \xh\p \xR $ & $\s 22$ & $28$\\
$[\s_1 \dd\s_1 \dd]$ & $\xh\p \xh\p \xh\p \xR $ & $\s 22$ & $39$\\
&&&\\

$(\s \d\s \d)$ & $\xo\p \xo\p \xo\p $ & $\s $ & $6$\\
$(\s \d\s \dd)$ & $\xo\p \xo\p \xh\p $ & $\s $ & $8$\\
$(\s \dd\s \dd)$ & $\xo\p \xh\p \xh\p $ & $\s $ & $7$\\
&&&\\

$(\bar{\s} \d\bar{\s} \d)$ & $\xo\m \xo\p \xo\p $ & $\s 22$ & $26$\\
$(\bar{\s} \d\bar{\s} \dd)$ & $\xo\m \xo\p \xh\p $ & $\s 22$ & $36$\\
$(\bar{\s} \dd\bar{\s} \dd)$ & $\xo\m \xh\p \xh\p $ & $\s 22$ & $29$\\
&&&\\

$(\s \d\s_0)$ & $\xo\p \xo\p \xo\m $ & $\s 22$ & $28$\\
$(\s \d\s_1)$ & $\xh\p \xo\p \xo\m $ & $\s 22$ & $40$\\
$(\s \dd\s_0)$ & $\xo\p \xh\p \xo\m $ & $\s 22$ & $32$\\
$(\s \dd\s_1)$ & $\xh\p \xh\p \xo\m $ & $\s 22$ & $41$\\
&&&\\

$(\bar{\s} \d\bar{\s}_0)$ & $\xo\m \xo\p \xo\m $ & $\s 22$ & $39$\\
$(\bar{\s} \d\bar{\s}_1)$ & $\xh\m \xo\p \xo\m $ & $\s 22$ & $46$\\
$(\bar{\s} \dd\bar{\s}_0)$ & $\xo\m \xh\p \xo\m $ & $\s 22$ & $45$\\
$(\bar{\s} \dd\bar{\s}_1)$ & $\xh\m \xh\p \xo\m $ & $\s 22$ & $41$\\
&&&\\

$(\s_0 \s_0)$ & $\xo\p \xo\m \xo\m $ & $22$ & $3$\\
$(\s_1 \s_1)$ & $\xh\p \xo\m \xo\m $ & $22$ & $5$\\
&&&\\

$(\bar{\s}_0 \bar{\s}_0)$ & $\xo\m \xo\m \xo\m $ & $\s 22$ & $27$\\
$(\bar{\s}_0 \bar{\s}_1)$ & $\xo\m \xo\m \xh\m $ & $\s 22$ & $37$\\
$(\bar{\s}_1 \bar{\s}_1)$ & $\xo\m \xh\m \xh\m $ & $\s 22$ & $30$\\
&&&\\
\end{supertabular}

\begin{supertabular}{|c|l|c|c|}

\subheader
{$\s\x$}
{$\s\up{P}\x\up{Z}$}
{$1=P^2=[P,Z^2]$}
{$P$ \, $Z$ \,\id}

$[\s \d\x_0]$ & $\xo\p \xo\p \xR $ & $\s 22$ & $38$\\
$[\s \d\x_1]$ & $\xo\p \xh\p \xR $ & $\s 22$ & $44$\\
$[\s \dd\x_0]$ & $\xh\p \xo\p \xR $ & $\s 22$ & $46$\\
$[\s \dd\x_1]$ & $\xh\p \xh\p \xR $ & $\s 22$ & $40$\\
&&&\\

$(\s \d\x)$ & $\xo\p \xo\p $ & $\s $ & $8$\\
$(\s \dd\x)$ & $\xh\p \xo\p $ & $\s $ & $9$\\
&&&\\

$(\s \d\bar{\x})$ & $\xo\p \xo\m $ & $\s 22$ & $31$\\
$(\s \dd\bar{\x})$ & $\xh\p \xo\m $ & $\s 22$ & $33$\\
&&&\\

$(\s_0 \x_0)$ & $\xo\m \xo\p $ & $\s 22$ & $30$\\
$(\s_0 \x_1)$ & $\xo\m \xh\p $ & $\s 22$ & $34$\\
$(\s_1 \x)$ & $\xo\m \xq\p $ & $\s 22$ & $43$\\
&&&\\

$(\s  \bar{\x})$ & $\xo\m \xo\m $ & $22$ & $5$\\
&&&\\

\subheader
{$\x\x$}
{$\x\up{Y}\x\up{Z}$}
{$1=Y^2 Z^2$}
{$Y$ {\hspace{3pt}} $Z$ {\hspace{3pt}} \id}

$[\x_0\x_0]$ & $\xo\p \xo\p \xR $ & $\s 22$ & $26$\\
$[\x_0\x_1]$ & $\xo\p \xh\p \xR $ & $\s 22$ & $36$\\
$[\x_1\x_1]$ & $\xh\p \xh\p \xR $ & $\s 22$ & $31$\\
&&&\\

$(\x\x_0)$ & $\xo\p \xo\p $ & $\s $ & $7$\\
$(\x\x_1)$ & $\xo\p \xh\p $ & $\s $ & $9$\\
&&&\\

$(\bar{\x} \x_0)$ & $\xo\m \xo\p $ & $\s 22$ & $29$\\
$(\bar{\x} \x_1)$ & $\xo\m \xh\p $ & $\s 22$ & $33$\\
&&&\\

$(\bar{\x}\bar{\x})$ & $\xo\m \xo\m $ & $22$ & $4$\\
&&&\\

\subheader
{$\o$}
{$\o\up{X}\up{Y}$}
{$1=[X,Y]$}
{$X$ {\hspace{3pt}} $Y$ {\hspace{3pt}} \id}

$[\o_0 ]$ & $\xo\p \xo\p \xR $ & $\s $ & $6$\\
$[\o_1 ]$ & $\xh\p \xh\p \xR $ & $\s $ & $8$\\
&&&\\

$(\o )$ & $\xo\p \xo\p $ & $1$ & $1$\\
&&&\\

$(\bar{\o} _0)$ & $\xo\m \xo\m $ & $\s $ & $7$\\
$(\bar{\o} _1)$ & $\xo\m \xh\m $ & $\s $ & $9$\\
&&&\\
\end{supertabular}

\end{center}

%
%
%
\twocolumn
\begin{center}

\tablehead{\hline\multicolumn{3}{|c|}{\em Table 2a continued}\\
\hline Primary & International & Secondary\\
name & no.\ and name&names \\\hline\hline}

\tablefirsthead{\hline\multicolumn{3}{|c|}{\bf Table 2a (\today)}\\
\hline Primary & International & Secondary\\
name & no.\ and name&names \\\hline\hline}

\tabletail{\hline}

\begin{supertabular}{|l|l|l|}
\multicolumn{3}{|l|}{Point Group $ \s 432$, \quad nos.\ 221-230}\\\hline
$8\up{\o}\dd2$ & $229.Im\b{3}m$ & $[\s\d4\d4\dd2]\dd3$, $[2_1\s\d2\d2]\dd6$\\
$8\up{\o}$ & $223.Pm\b{3}n$ & $[\s\d4\dd4\d2]\dd3$, $[\s\d2\d2\d2\d2]\dd6$\\
$8\up{\o\o} $ & $222.Pn\b{3}n$ & $(\s4_04\dd2)\dd3$, $(2\bar{\s}_12_02_0)\dd6$\\ 
$4\up{\m}\dd2$ & $221.Pm\b{3}m$ & $[\s\d4\d4\d2]\dd3$, $[\s\d2\d2\d2\d2]\dd6$\\ 
$4\up{\p}\dd2$ & $224.Pn\b{3}m$ & $(\s4_24\d2)\dd3$, $(2\bar{\s}_12_02_0)\dd6$\\
$4\up{\m\m}$ & $226.Fm\b{3}c$ & $[\s\d4\dd4\dd2]\dd3$, $[\s\d2\d2\dd2\dd2]\dd6$\\
$4\up{\p\p}$ & $228.Fd\b{3}c$ & $(\s4_14\dd2)\dd3$, $(2\bar{\s}2_02_1)\dd6$\\ 
$2\up{\m}\dd2$ & $225.Fm\b{3}m$ & $[\s\d4\d4\dd2]\dd3$, $[\s\d2\d2\dd2\dd2]\dd6$\\
$2\up{\p}\dd2$ & $227.Fd\b{3}m$ & $(\s4_14\d2)\dd3$, $(2\bar{\s}2_02_1)\dd6$\\ 
$8\up{\o}/4$ & $230.Ia\b{3}d$ & $(\s4_14\dd2)\dd3$, $(\s2_12\dd2\dd2)\dd6$\\ 
&&\\
\hline \multicolumn{3}{|l|}{Point Group $ 432$, \quad nos.\ 207-214}\\
\hline
$8\up{\p \o} $ & $211.I432$ &$(\s4_24_02_1)\dd3$, $(2_1\s2_02_0)\dd6$\\ 
$4\up{\p} $ & $208.P4_232$ & $(\s4_24_22_0)\dd3$, $(\s2_02_02_02_0)\dd6$\\
$4\up{\o \m} $ & $207.P432$ &$(\s4_04_02_0)\dd3$, $(\s2_02_02_02_0)\dd6$\\
$2\up{\p} $ & $210.F4_132$ & $(\s4_34_12_0)\dd3$, $(\s2_02_12_02_1)\dd6$\\
$2\up{\o \m} $ & $209.F432$ &$(\s4_24_02_1)\dd3$, $(\s2_02_12_02_1)\dd6$\\
$4\up{\p}/4$ & $214.I4_132$ & $(\s4_34_12_0)\dd3$, $(2_0\s2_12_1)\dd6$\\
$2\up{\p}/4$ & $\ddagger$ $212.P4_332$ & $(4_1\s2_1)\dd3$, $(2_12_1\bar{\x})\dd6$\\
&&\\
\hline \multicolumn{3}{|l|}{Point Group $ 3\s 2$, \quad nos.\ 200-206}\\
\hline
$8\up{\m \o} $ & $204.Im\b{3}$ &$[2_1\s\d2\d2]\dd3$\\
$4\up{\m} $ & $200.Pm\b{3}$ & $[\s\d2\d2\d2\d2]\dd3$\\
$4\up{\o \p} $ & $201.Pn\b{3}$ & $(2\bar{\s}_12_02_0)\dd3$\\
$2\up{\m} $ & $202.Fm\b{3}$ & $[\s\d2\d2\dd2\dd2]\dd3$\\
$2\up{\o \p} $ & $203.Fd\b{3}$ & $(2\bar{\s}2_02_1)\dd3$\\
$4\up{\m} /4$ & $206.Ia\b{3}$ & $(\s2_12\dd2\dd2)\dd3$\\
$2\up{\m} /4 $ & $205.Pa\b{3}$ & $(2_12\bar{\s}\dd)\dd3$\\
&&\\
\hline \multicolumn{3}{|l|}{Point Group $ \s 332$, \quad nos.\ 215-220}\\
\hline
$4\up{\o} \dd2$ & $217.I\b{4}3m$ &$(\s.4\n4\dd2)\dd3$, $(2_1\s2_02_0)\dd6$\\ 
$4\up{\o} $ & $218.P\b{4}3n$ & $(\s4\dd4\n2_0)\dd3$, $(\s2_02_02_02_0)\dd6$\\
$2\up{\o} \dd2$ & $215.P\b{4}3m$ &$(\s4\d4\n2_0)\dd3$, $(\s2_02_02_02_0)\dd6$\\
$2\up{\o \o} $ & $219.F\b{4}3c$ & $(\s4\dd4\n2_1)\dd3$, $(\s2_02_12_02_1)\dd6$\\
$1\up{\o} \dd2$ & $216.F\b{4}3m$ &$(\s4\d4\n2_1)\dd3$, $(\s2_02_12_02_1)\dd6$\\
$4\up{\o} /4$ & $220.I\b{4}3d$ & $(4\bar{\s}2_1)\dd3$, $(2_0\s2_12_1)\dd6$\\ 
&&\\
\hline \multicolumn{3}{|l|}{Point Group $ 332$, \quad nos.\ 195-199}\\
\hline
$4\up{\o\o} $ & $197.I23$ & $(2_1\s2_02_0)\dd3$\\
$2\up{\o} $ & $195.P23$ & $(\s2_02_02_02_0)\dd3$\\
$1\up{\o} $ & $196.F23$ & $(\s2_02_12_02_1)\dd3$\\
$2\up{\o} /4$ & $199.I2_13$ &$(2_0\s2_12_1)\dd3$\\
$1\up{\o} /4$ & $198.P2_13$ &$(2_12_1\bar{\x})\dd3$\\
&&\\
\end{supertabular}
{\tt Memo: (1) Compare this carefully with Olaf's file and Table2b.
(2) different groups have the same :6 secondary name, so is
the latter really a name?}

\end{center}

%
%
%
%
%
%
\twocolumn
\begin{center}

\tablehead{\hline\multicolumn{3}{|c|}{\it Table~2b continued}\\
\hline Primary & International &  Secondary\\
name & no.\ and name&names \\\hline\hline}

\tablefirsthead{\hline\multicolumn{3}{|c|}{\bf Table~2b (\today)}\\
\hline Primary & International &  Secondary\\
name & no.\ and name&names \\\hline\hline}

\tabletail{\hline}

\begin{supertabular}{|l|l|l|}
\hline \multicolumn{3}{|l|}{Point Group $ \s 226$, \quad nos.\ 191-194}\\ \hline
$[\s \d6\d3\d2] $ & $191.P6/mmm$ &\\
$[\s \dd6\d3\d2] $ & $194.P6_3/mmc$ &\\
$[\s \d6\dd3\dd2] $ & $193.P6_3/mcm$ &\\
$[\s \dd6\dd3\dd2] $ & $192.P6/mcc$ &\\
&&\\
\hline \multicolumn{3}{|l|}{Point Group $ 226$, \quad nos.\ 177-182}\\ \hline
$(\s 6_0 3_0 2_0) $ & $177.P622$ &\\
$(\s 6_1 3_1 2_1) $ & $\ddagger$ $178.P6_122$ &\\
$(\s 6_2 3_2 2_0) $ & $\ddagger$ $180.P6_222$ &\\
$(\s 6_3 3_0 2_1) $ & $182.P6_322$ &\\
&&\\
\hline \multicolumn{3}{|l|}{Point Group $ \s 66$, \quad nos.\ 183-186}\\ \hline
$(\s \d6\d3\d2) $ & $183.P6mm$ &\\
$(\s \dd6\d3\d2) $ & $186.P63mc$ &\\
$(\s \d6\dd3\dd2) $ & $185.P6_3cm$ &\\
$(\s \dd6\dd3\dd2) $ & $184.P6cc$ &\\
&&\\
\hline \multicolumn{3}{|l|}{Point Group $ 6\s $, \quad nos.\ 175-176}\\ \hline
$[6_0 3_0 2_0] $ & $175.P6/m$ &\\
$[6_3 3_0 2_1] $ & $176.P6_3/m$ &\\
&&\\
\hline \multicolumn{3}{|l|}{Point Group $ 66$, \quad nos.\ 168-173}\\ \hline
$(6_0 3_0 2_0) $ & $168.P6$ &\\
$(6_1 3_1 2_1) $ & $\ddagger$ $169.P6_1$ &\\
$(6_2 3_2 2_0) $ & $\ddagger$ $171.P6_2$ &\\
$(6_3 3_0 2_1) $ & $173.P6_3$ &\\
&&\\
\hline \multicolumn{3}{|l|}{Point Group $ 2\s 3$, \quad nos.\ 162-167}\\ \hline
$(\s 6\d3\d2) $ & $164.P\b{3}m1$ &\\
$(\s 6\dd3\dd2) $ & $165.P\b{3}c1$ &\\
&&\\
$(\s \d6\n 3_0 2) $ & $162.P\b{3}1m$ &\\
$(\s \dd6\n 3_0 2) $ & $163.P\b{3}1c$ &\\
$(\s \d6\n 3_1 2) $ & $166.R\b{3}m$ &\\
$(\s \dd6\n 3_1 2) $ & $167.R\b{3}c$ &\\
&&\\
\hline \multicolumn{3}{|l|}{Point Group $ 3\x$, \quad nos.\ 147-148}\\ \hline
$(6\n 3_0 2) $ & $147.P\b{3}$ &\\
$(6\n 3_1 2) $ & $148.R\b{3}$ &\\
&&\\
\hline \multicolumn{3}{|l|}{Point Group $ \s 224$, \quad nos.\ 123-142}\\ \hline
$[\s \d4\d4\d2] $ & $123.P4/mmm$ &\\
$[\s \d4\d4\dd2] $ & $139.I4/mmm$ &\\
$[\s \d4\dd4\d2] $ & $131.P4_2/mmc$ &\\
$[\s \d4\dd4\dd2] $ & $140.I4/mcm$ &\\
$[\s \dd4\d4\dd2] $ & $132.P4_2/mcm$ &\\
$[\s \dd4\dd4\dd2] $ & $124.P4/mcc$ &\\
&&\\
$(\s 4\d4\d2) $ & $129.P4/nmm$ &\\
$(\s 4\d4\dd2) $ & $137.P4_2/nmc$ &\\
$(\s 4\dd4\d2) $ & $138.P4_2/ncm$ &\\
$(\s 4\dd4\dd2) $ & $130.P4/ncc$ &\\
&&\\
$(\s 4_0 4\d2) $ & $125.P4/nbm$ &\\
$(\s 4_0 4\dd2) $ & $126.P4/nnc$ &\\
$(\s 4_1 4\d2) $ & $141.I4_1/amd$ &\\
$(\s 4_1 4\dd2) $ & $142.I4_1/acd$ &\\
$(\s 4_2 4\d2) $ & $134.P4_2/nnm$ &\\
$(\s 4_2 4\dd2) $ & $133.P4_2/nbc$ &\\
&&\\
$[4_0 \s \d2] $ & $127.P4/mbm$ &\\
$[4_0 \s \dd2] $ & $128.P4/mnc$ &\\
$[4_2 \s \d2] $ & $136.P4_2/mnm$ &\\
$[4_2 \s \dd2] $ & $135.P4_2/mbc$ &\\
&&\\
\hline \multicolumn{3}{|l|}{Point Group $ 224$,  \quad nos.\ 89-98}\\ \hline
$(\s 4_0 4_0 2_0) $ & $89.P422$ &\\
$(\s 4_1 4_1 2_1) $ & $\ddagger$ $91.P4_122$ &\\
$(\s 4_2 4_2 2_0) $ & $93.P4_222$ &\\
$(\s 4_2 4_0 2_1) $ & $97.I422$ &\\
$(\s 4_3 4_1 2_0) $ & $98.I4_122$ &\\
&&\\
$(4_0 \s 2_0) $ & $90.P42_12$ &\\
$(4_1 \s 2_1) $ & $\ddagger$ $92.P4_12_12$ &\\
$(4_2 \s 2_0) $ & $94.P4_22_12$ &\\
&&\\
\hline \multicolumn{3}{|l|}{Point Group $ \s 44$, \quad nos.\ 99-110}\\ \hline
$(\s \d4\d4\d2) $ & $99.P4mm$ &\\
$(\s \d4\d4\dd2) $ & $107.I4mm$ &\\
$(\s \d4\dd4\d2) $ & $105.P4_2mc$ &\\
$(\s \d4\dd4\dd2) $ & $108.I4cm$ &\\
$(\s \dd4\d4\dd2) $ & $101.P4_2cm$ &\\
$(\s \dd4\dd4\dd2) $ & $103.P4cc$ &\\
&&\\
$(4_0 \s \d2) $ & $100.P4bm$ &\\
$(4_0 \s \dd2) $ & $104.P4nc$ &\\
$(4_1 \s \d2) $ & $109.I4_1md$ &\\
$(4_1 \s \dd2) $ & $110.I41cd$ &\\
$(4_2 \s \d2) $ & $102.P4_2nm$ &\\
$(4_2 \s \dd2) $ & $106.P4_2bc$ &\\
&&\\
\end{supertabular}

\tablefirsthead{\hline\multicolumn{3}{|c|}{\it Table~2b continued}\\
\hline Primary & International &  Secondary\\
name & no.\ and name&names \\\hline\hline}

\begin{supertabular}{|l|l|l|}
\hline \multicolumn{3}{|l|}{Point Group $ 4\s $, \quad nos.\ 83-88}\\ \hline
$[4_0 4_0 2_0] $ & $83.P4/m$ &\\
$[4_2 4_2 2_0] $ & $84.P4_2/m$ &\\
$[4_2 4_0 2_1] $ & $87.I4/m$ &\\
&&\\
$(4\n 4_0 2) $ & $85.P4/n$ &\\
$(4\n 4_1 2) $ & $88.I4_1/a$ &\\
$(4\n 4_2 2) $ & $86.P4_2/n$ &\\
&&\\
\hline \multicolumn{3}{|l|}{Point Group $ 44$, \quad nos.\ 75-80}\\ \hline
$(4_0 4_0 2_0) $ & $75.P4$ &\\
$(4_1 4_1 2_1) $ & $\ddagger$ $76.P4_1$ &\\
$(4_2 4_2 2_0) $ & $77.P4_2$ &\\
$(4_2 4_0 2_1) $ & $79.I4$ &\\
$(4_3 4_1 2_0) $ & $80.I4_1$ &\\
&&\\
\hline \multicolumn{3}{|l|}{Point Group $ 2\s 2$, \quad nos.\ 111-122}\\ \hline
$(\s \d4\n 4\d2) $ & $115.P\b{4}m2$ &\\
$(\s .4\n 4\dd2) $ & $121.I\b{4}2m$ &\\
$(\s \dd4\n 4\dd2) $ & $116.P\b{4}c2$ &\\
&&\\
$(\s 4\d4\n 2_0) $ & $111.P\b{4}2m$ &\\
$(\s 4\dd4\n 2_0) $ & $112.P\b{4}2c$ &\\
$(\s 4\d4\n 2_1) $ & $119.I\b{4}m2$ &\\
$(\s 4\dd4\n 2_1) $ & $120.I\b{4}c2$ &\\
&&\\
$(4 \bar{\s} \d2) $ & $113.P\b{4}2_1m$ &\\
$(4 \bar{\s} \dd2) $ & $114.P\b{4}2_1c$ &\\
&&\\
$(4 \bar{\s}_0 2_0) $ & $117.P\b{4}b2$ &\\
$(4 \bar{\s}_1 2_0) $ & $118.P\b{4}n2$ &\\
$(4 \bar{\s} 2_1) $ & $122.I\b{4}2d$ &\\
&&\\
\hline \multicolumn{3}{|l|}{Point Group $ 2\x$, \quad nos.\ 81-82}\\ \hline
$(4\n 4\n 2_0) $ & $81.P\b{4}$ &\\
$(4\n 4\n 2_1) $ & $82.I\b{4}$ &\\
&&\\
\hline \multicolumn{3}{|l|}{Point Group $ \s 223$, \quad nos.\ 187-190}\\ \hline
$[\s \d3\d3\d3] $ & $187.P\b{6}m2$ &\\
$[\s \dd3\dd3\dd3] $ & $188.P\b{6}c2$ &\\
&&\\
$[3_0 \s \d3] $ & $189.P\b{6}2m$ &\\
$[3_0 \s \dd3] $ & $190.P\b{6}2c$ &\\
&&\\
\hline \multicolumn{3}{|l|}{Point Group $ 223$, \quad nos.\ 149-155}\\ \hline
$(\s 3_0 3_0 3_0) $ & $149.P312$ &\\
$(\s 3_1 3_1 3_1) $ & $\ddagger$ $151.P3_112$ &\\
$(\s 3_0 3_1 3_2) $ & $155.R32$ &\\
&&\\
$(3_0 \s 3_0) $ & $150.P321$ &\\
$(3_1 \s 3_1) $ & $\ddagger$ $152.P3_121$ &\\
&&\\
\end{supertabular}

\begin{supertabular}{|l|l|l|}
\hline \multicolumn{3}{|l|}{Point Group $ \s 33$, \quad nos.\ 156-161}\\ \hline
$(\s \d3\d3\d3) $ & $156.P3m1$ &\\
$(\s \dd3\dd3\dd3) $ & $158.P3c1$ &\\
&&\\
$(3_0 \s \d3) $ & $157.P31m$ &\\
$(3_0 \s \dd3) $ & $159.P31c$ &\\
$(3_1 \s \d3) $ & $160.R3m$ &\\
$(3_1 \s \dd3) $ & $161.R3c$ &\\
&&\\
\hline \multicolumn{3}{|l|}{Point Group $ 3\s $, \quad no.\ 174}\\ \hline
$[3_0 3_0 3_0] $ & $174.P\b{6}$ &\\
&&\\
\hline \multicolumn{3}{|l|}{Point Group $ 33$, \quad nos.\ 143-146}\\ \hline
$(3_0 3_0 3_0) $ & $143.P3$ &\\
$(3_0 3_1 3_2) $ & $146.R3$ &\\
$(3_1 3_1 3_1) $ & $\ddagger$ $144.P3_1$ &\\
&&\\
\hline \multicolumn{3}{|l|}{Point Group $ \s 222$, \quad nos.\ 47-74}\\ \hline
$[\s \d2\d2\d2\d2] $ & $47.Pmmm$ &\\
$[\s \d2\d2\dd2\dd2] $ & $69.Fmmm$ &\\
$[\s \dd2\dd2\dd2\dd2] $ & $49.Pccm$ & $(\s 2_0 2_0 2\d2)$\\
&&\\
$(\s 2_0 2\d2\d2) $ & $67.Cmma$ & $[\s \d2\dd2\dd2\dd2]$\\
$(\s 2_0 2\dd2\dd2) $ & $68.Ccca$ & $(\s 2_0 2_1 2\dd2)$\\
$(\s 2_1 2\d2\d2) $ & $74.Imma$ & $[2_0 \s \d2\dd2]$\\
$(\s 2_1 2\dd2\dd2) $ & $73.Ibca$ &\\
&&\\
$[2_0 \s \d2\d2] $ & $65.Cmmm$ & $[\s \d2\d2\d2\dd2]$\\
$[2_0 \s \dd2\dd2] $ & $66.Cccm$ & $(\s 2_0 2_1 2\d2)$\\
$[2_1 \s \d2\d2] $ & $71.Immm$ &\\
$[2_1 \s \dd2\dd2] $ & $72.Ibam$ & $(\s 2_0 2\d2\dd2)$\\
&&\\
$(2 \bar{\s} \d2\d2)$ & $59.Pmmn$ & $[2_1 2_1 \s \d]$\\
$(2 \bar{\s} \dd2\dd2)$ & $56.Pccn$ & $(2_1 2 \bar{\s}_0)$\\
&&\\
$(2 \bar{\s}_0 2_0 2_0)$ & $50.Pban$ & $(\s 2_0 2_0 2\dd2)$\\
$(2 \bar{\s}_1 2_0 2_0) $ & $48.Pnnn$ &\\
$(2 \bar{\s} 2_0 2_1) $ & $70.Fddd$ &\\
&&\\
$[2_0 2_0 \x_0] $ & $55.Pbam$ & $(\s 2\d2\dd2\d2)$\\
$[2_0 2_0 \x_1] $ & $58.Pnnm$ & $(2_1 \s 2\d2)$\\
&&\\
$[2_0 2_0 \s \d] $ & $51.Pmma$ & $[\s \d2\dd2\d2\dd2], (\s 2\d2\d2\d2)$\\
$[2_0 2_0 \s \dd] $ & $53.Pmna$ & $(\s 2_1 2_1 2\d2), (2_0 \s 2\d2)$\\
$[2_0 2_1 \s \d] $ & $63.Cmcm$ & $(\s 2\d2\d2\dd2), [2_1 \s \d2\dd2]$\\
$[2_0 2_1 \s \dd] $ & $64.Cmca$ & $(\s 2\d2\dd2\dd2), (\s 2_1 2\d2\dd2)$\\
&&\\
$(2_0 2 \bar{\s} \d) $ & $57.Pbcm$ & $(\s 2\dd2\d2\dd2), [2_1 2_1 \s \dd]$\\
$(2_0 2 \bar{\s} \dd) $ & $60.Pbcn$ & $(2_1 \s 2\dd2), (2_1 2 \bar{\s}_1)$\\
$(2_1 2 \bar{\s} \d) $ & $62.Pnma$ & $(2 \bar{\s} \d2\dd2), [2_1 2_1 \x]$\\
$(2_1 2 \bar{\s} \dd) $ & $61.Pbca$ &\\
&&\\
$(2_0 2 \bar{\s}_0) $ & $54.Pcca$ & $(\s 2\dd2\dd2\dd2), (\s 2_1 2_1 2\dd2)$\\
$(2_0 2 \bar{\s}_1) $ & $52.Pnna$ & $(2_0 \s 2\dd2), (2 \bar{\s} 2_1 2_1)$\\
&&\\
\end{supertabular}

\begin{supertabular}{|l|l|l|}
\hline \multicolumn{3}{|l|}{Point Group $ 222$, \quad nos.\ 16-24}\\ \hline
$(\s 2_0 2_0 2_0 2_0)$ & $16.P222$ &\\
$(\s 2_0 2_1 2_0 2_1)$ & $22.F222$ &\\
$(\s 2_1 2_1 2_1 2_1)$ & $17.P222_1$ & $(2_0 2_0 \s )$\\
&&\\
$(2_0 \s 2_0 2_0) $ & $21.C222$ & $(\s 2_0 2_0 2_1 2_1)$\\
$(2_1 \s 2_0 2_0) $ & $23.I222$ &\\
$(2_0 \s 2_1 2_1) $ & $24.I2_12_12_1$ &\\
$(2_1 \s 2_1 2_1) $ & $20.C222_1$ & $(2_0 2_1 \s )$\\
&&\\
$(2_0 2_0 \bar{\x}) $ & $18.P2_12_12$ & $(2_1 2_1 \s )$\\
$(2_1 2_1 \bar{\x}) $ & $19.P2_12_12_1$ &\\
&&\\
\hline \multicolumn{3}{|l|}{Point Group $ \s 22$, \quad nos.\ 25-46}\\ \hline
$(\s \d2\d2\d2\d2) $ & $25.Pmm2$ & $[\s_0 \d\s_0 \d]$\\
$(\s \d2\d2\d2\dd2) $ & $38.Amm2$ & $[\s_1 \d\s_1 \d], [\s \d\x_0]$\\
$(\s \d2\d2\dd2\dd2) $ & $42.Fmm2$ & $[\s_1 \d\s_1 \dd]$\\
$(\s \d2\dd2\d2\dd2) $ & $26.Pmc2_1$ & $(\bar{\s} \d\bar{\s} \d), [\x_0\x_0]$\\
$(\s \d2\dd2\dd2\dd2) $ & $39.Abm2$ & $[\s_1 \dd\s_1 \dd],
(\bar{\s} \d\bar{\s}_0)$\\
$(\s \dd2\dd2\dd2\dd2) $ & $27.Pcc2$ & $(\bar{\s}_0 \bar{\s}_0)$\\
&&\\
$(2_0 \s \d2\d2) $ & $35.Cmm2$ & $[\s_0 \d\s_0 \dd]$\\
$(2_0 \s \d2\dd2) $ & $46.Ima2$ & $(\bar{\s} \d\bar{\s}_1), [\s \dd\x_0]$\\
$(2_0 \s \dd2\dd2) $ & $37.Ccc2$ & $(\bar{\s}_0 \bar{\s}_1)$\\
$(2_1 \s \d2\d2) $ & $44.Imm2$ & $[\s \d\x_1]$\\
$(2_1 \s \d2\dd2) $ & $36.Cmc2_1$ & $(\bar{\s} \d\bar{\s} \dd), [\x_1\x_1]$\\
$(2_1 \s \dd2\dd2) $ & $45.Iba2$ & $(\bar{\s} \dd\bar{\s}_0)$\\
&&\\
$(2_0 2_0 \s \d) $ & $28.Pma2$ & $[\s_0 \dd\s_0 \dd], (\s \d\s_0)$\\
$(2_0 2_0 \s \dd) $ & $30.Pnc2$ & $(\bar{\s}_1 \bar{\s}_1), (\s_0 \x_0)$\\
$(2_0 2_1 \s \d) $ & $40.Ama2$ & $(\s \d\s_1), [\s \dd\x_1]$\\
$(2_0 2_1 \s \dd) $ & $41.Aba2$ & $(\s \dd\s_1), (\bar{\s} \dd\bar{\s}_1)$\\
$(2_1 2_1 \s \d) $ & $31.Pmn2_1$ & $(\s \d\bar{\x}), [\x_0\x_1]$\\
$(2_1 2_1 \s \dd) $ & $29.Pca2_1$ & $(\bar{\s} \dd\bar{\s} \dd),
(\bar{\x} \x_0)$\\
&&\\
$(2_0 2_0 \x_0) $ & $32.Pba2$ & $(\s \dd\s_0)$\\
$(2_0 2_0 \x_1) $ & $34.Pnn2$ & $(\s_0 \x_1)$\\
$(2_0 2_1 \x) $ & $43.Fdd2$ & $(\s_1 \x)$\\
$(2_1 2_1 \x) $ & $33.Pa2_1$ & $(\s \dd\bar{\x}), (\bar{\x} \x_1)$\\
&&\\
\hline \multicolumn{3}{|l|}{Point Group $ 2\s $, \quad nos.\ 10-15}\\ \hline
$[2_0 2_0 2_0 2_0] $ & $10.P2/m$ & $(\s 2\d2\n 2\d2)$\\
$[2_0 2_0 2_1 2_1] $ & $12.C2/m$ & $(\s 2\d2\n 2\dd2), (2 \bar{\s} 2\d2)$\\
$[2_1 2_1 2_1 2_1] $ & $11.P2_1/m$ & $(2\n 2 \s \d)$\\
&&\\
$(2_0 2_0 2\n 2) $ & $13.P2/c$ & $(\s 2\dd2\n 2\dd2), (2\n 2 \s_0)$\\
$(2_0 2_1 2\n 2) $ & $15.C2/c$ & $(2 \bar{\s} 2\dd2), (2\n 2 \s_1)$\\
$(2_1 2_1 2\n 2) $ & $14.P2_1/c$ & $(2\n 2 \s \dd), (2\n 2 \x)$\\
&&\\
\hline \multicolumn{3}{|l|}{Point Group $ 22$, \quad nos.\ 3-5}\\ \hline
$(2_0 2_0 2_0 2_0) $ & $3.P2$ & $(\s_0 \s_0)$\\
$(2_0 2_0 2_1 2_1) $ & $5.C2$ & $(\s_1 \s_1), (\s \bar{\x})$\\
$(2_1 2_1 2_1 2_1) $ & $4.P2_1$ & $(\bar{\x}\bar{\x})$\\
&&\\
\end{supertabular}

\begin{supertabular}{|l|l|l|}
\hline \multicolumn{3}{|l|}{Point Group $ \x$, \quad no.\ 2}\\ \hline
$(2\n 2\n 2\n 2) $ & $2.P\b{1}$ &\\
&&\\
\hline \multicolumn{3}{|l|}{Point Group $ \s $, \quad nos.\ 6-9}\\ \hline
$[\o_0 ] $ & $6.Pm$ & $(\s \d\s \d)$\\
$[\o_1 ] $ & $8.Cm$ & $(\s \d\s \dd), (\s \d\x)$\\
&&\\
$(\bar{\o}_0) $ & $7.Pc$ & $(\s \dd\s \dd),
(\x\x_0)$\\
$(\bar{\o}_1) $ & $9.Cc$ & $(\s \dd\x), (\x\x_1)$\\
&&\\
\hline \multicolumn{3}{|l|}{Point Group $ 1$, \quad no.\ 1}\\ \hline
$(\o ) $ & $1.P1$ &\\
&&\\
\end{supertabular}

\end{center}


\end{document}